\documentclass[a4paper]{scrartcl}
\usepackage{amsmath}
\usepackage{amssymb}
\usepackage{amsthm}
\usepackage[utf8]{inputenc}
\usepackage{geometry}
\usepackage{hyperref}
\usepackage[nocites,norefs]{refcheck}
\usepackage[usenames,dvipsnames]{xcolor}

\DeclareMathOperator{\Sp}{Sp}
\DeclareMathOperator{\Ep}{Ep}
\DeclareMathOperator{\GL}{GL}
\DeclareMathOperator{\Up}{Up}
\DeclareMathOperator{\U}{U}
\DeclareMathOperator{\PC}{PC}
\DeclareMathOperator{\Ann}{Ann}
\DeclareMathOperator{\Center}{Center}
\DeclareMathOperator{\height}{ht}
\DeclareMathOperator{\rk}{rk}
\DeclareMathOperator{\Int}{Int}
\DeclareMathOperator{\Aut}{Aut}
\DeclareMathOperator{\es}{es}

\newcommand{\red}[1]{\textcolor{red}{#1}}
\renewcommand{\red}[1]{}

\renewcommand{\tilde}[1]{\widetilde{#1}}
\renewcommand{\bar}[1]{\overline{#1}}

\newtheorem{lemma}{Lemma}
\newtheorem{corollary}[lemma]{Corollary}
\newtheorem{theorem}[lemma]{Theorem}

\author{Alexander Shchegolev\footnote{Research is supported by the Russian Science Foundation grant No17-11-01261.}}
\title{Bijective PC-maps of the unipotent radical of the Borel subgroup of the
classical symplectic group}
\date{Saint Petersburg State University}

\begin{document}

\maketitle

\begin{abstract}
	Every commutator preserving bijection of the unipotent radical $\Up(2n, R)$ of the Borel subgroup of the classical
	symplectic group of rank at least 4 over a field $F$ such that $6F=F$ is shown to be the composition
	of a standard automorphism of $\Up(2n,R)$ and a central map. The latter is a bijection which acts as the right
	multiplication by a matrix from the center of $\Up(2n,F)$.
\end{abstract}


Linear preserver problems (LPP) is a vast class of problems in linear algebra and related areas.
The general shape of such a problem is to describe all maps between linear spaces, groups,
Lie algebras etc preserving certain properties, operations, subspaces etc. The class of
the maps one wishes to describe is usually narrowed by some natural condition. Most of the 
original papers on LPP's dealt only with linear operators, while the others focused on
bijective but not necessarily linear maps. This research is focused on bijective commutators
preserving maps. There is a number of reviews on the state of the area of LPP's
as well as general methods thereof, e.g. \cite{LiPierceLPP2001, PlatDokLinPres95, GutLiSemLPP2000}.
Here we limit ourselves to mentioning several papers which are most closely related to this one 
both in terms of the specific problem considered and methods used. A description of bijective
commutators preserving maps on the subgroup of upper triangular matrices in $\GL(2n, F)$, where
$F\neq\mathbb{F}_3$ and $\operatorname{char}{F}\neq 2$, was given by Dengyin Wang, Shikun Ou and Wei Zhang in 
\cite{WangOuZhangPCB2010} up to so called almost identity maps. This classification was later 
carried over to the case of upper unitriangular matrices by Dengyin Wang, Huixiang Zhai and 
Meilin Chen in \cite{WangZhaiChenPCU2011} under the usual assumption that $\operatorname{char}{F}\neq 2$. 
Further, in \cite{SlowPCInf2013} Roksana Slowik gave a similar description for $\PC$-maps on infinite upper triangular
and unitriangular subgroups of $\GL(\infty,F)$, again up to almost identity maps. It was unclear
what almost identity maps actually look like up until the recent article \cite{StepHolPC2015} by Alexei Stepanov and
Waldemar Holubowski which shows that these are precisely central maps, i.e. maps acting as
multiplication on the right by an element of the center of the group. Group theoretic
properties of the group of $\PC$-maps are also studied there. Next, in \cite{OuPCBSp2012}
Shikun Ou describes $\PC$-maps of the Borel subgroup of the classical symplectic group by reducing
this problem to the linear case and the description of automorphisms of the group in consideration
given by Ou and Wnag in \cite{OuWangAutomBSp2008}. A closely related series of results concerning 
(not necessarily linear) zero Lie product preservers on Borel subalgebras of simple Lie algebras 
was obtained by Dengyin Wang in \cite{WangLiezero2014}.

The central result of this paper is the following theorem. The standard $\PC$-maps
mentioned in this theorem are defined in Section \ref{sec:standard}.
\begin{theorem}
	\label{th:main}
	Let $n \geq 4$ and $F$ be a field such that $6F=F$. 
	Let $\Up$ denote the unipotent radical of the Borel subgroup
	of the classical symplectic group $\Sp(2n,F)$ of rank $n$ over $F$. Then any
	commutator preserving bijective map $\phi$ is a composition
	of an inner automorphism, an extremal automorphism, a diagonal automorphism, a semi-diagonal automorphism,
	a field automorphism and a central $\PC$-map.
\end{theorem}

\section{Preliminaries}

\paragraph{Symplectic group.}
\label{sec:prelim}

Fix a natural number $n$ and a field $F$. Let $\Sp(2n, R)$ denote the classical symplectic group
of rank $n$ over $F$. Let $\Up=\Up(2n,F)$ denote the subgroup of upper unitriangular matrices in $\Sp(2n,F)$.
In terms of algebraic groups, $\Sp(2n, F)$ is the group of points of the adjoint Chevalley group 
of type $C_n$ over $F$ and $\Up(2n, R)$ the unipotent radical of the Borel subgroup of $\Sp(2n,F)$. 

We will use the ordinary set $\{1,\dots,2n\}$ to index the entries of matrices of $\Sp(2n, R)$ rather than 
the set $\{1,\dots,n,-n,\dots,-1\}$. The latter is used by many authors to simplify computations dealing with
elementary matrices. However in this exposition we use roots of the root system $C_n$ to index the elementary root 
unipotents while working with generators of $\Up(2n,F)$ and use explicitly written matrices
in direct computations which do not deal with root unipotents directly. The former do not involve matrices at all,
while the latter are performed easier when the linear ordering on the indexing set is straightforward.

\paragraph{Elementary root unipotents and the root system of type $C_n$.}

When it comes to calculating in classical linear groups, there has been a long going rivalry.
One perspective is to treat all matrices solely as matrices and compute all products, conjugates and commutators
explicitly as products, conjugates and commutators of matrices. The benefits of this approach include
the ease of verification of computation and a low entry level for a reader. The downside is obscurity.
Once you read any construction like ``\texttt{let $i,j,k,l,t,p,q$ be indices such that $i \neq p, 2n-q, l+k,
j \neq p, q, k, l+1$ and any triple of these indices contains an element not divisible by two; 
then the $(l+k+j-t, p+q-1)$'th
entry of the comutator $[[[a,b],c],d]$}'' it's almost impossible to keep all these conditions in mind.
Usually, such conditions simply mean that we consider some subsystem of a root system. Which leads us to
the second approach: working with elementary subgroups of linear groups in terms of generators (elementary root
unipotents; elementary semisimple elements) and relations (Chevalley commutator formulae/Steinberg relations).
This way of thinking is remarkably more suggestive. Once one understands the root system of the group in 
consideration, it gets remarkably easy to compute products and commutators of ``not too much'' of generators.
In particular, this approach works effectively for matrices contained in certain parabolic subgroups.
However when dealing with arbitrary matrices it's often easier to compute a certain entry of a commutator
explicitly, rather than in terms of generators. Long story short, in this paper we are trying to get the best
of two worlds and make use of both styles of computations. However, deep understanding of neither root systems
nor algebraic groups is required to read this paper. To facilitate this we provide a quick and dirty 
introduction for translating from the [subset of the] language of Chevalley groups used in this paper
to that of classical linear groups.

\textit{The root system $\bar{\Phi}$ of type $C_n$} is the following set of vectors in $n$-dimensional real vector space
\[
	\bar{\Phi} = \left\lbrace \pm 2e_i, \pm e_i \pm e_j \mid 1 \leq i \neq j \leq n \right\rbrace,
\]
where $e_i$ constitute the standard base of $\mathbb{R}^n$.
The sum of roots simply means the sum of vectors. We will call the roots 
\[
	\alpha_i = \begin{cases} e_i-e_{i+1} & 1 \leq i \leq n-1 \\ 2e_n & i = n \end{cases}
\]
\textit{simple roots}, the set of simple roots is denoted by $\Pi$. Every root $\alpha \in \bar{\Phi}$ can be uniquely decomposed
as integral linear combination of simple root with all coefficients being either nonnegative
or nonpositive. In the former case the root is called \textit{positive} and in the latter
\textit{negative}. We will denote the set of positive roots by $\Phi$, which is not
traditional, as $\Phi$ usually denotes the whole root system, while the set of positive roots is 
denoted by $\Phi^+$. However in this paper we will never
use all roots, but only positive ones. The roots of Euclidean length 2 are called long, while 
those of length $\sqrt{2}$ \textit{short}. The sets of long and short positive roots will be denoted
by $\Phi^L$ and $\Phi^S$, respectively. Summing up,
\begin{align*}
	\Phi &= \left\lbrace 2e_i, e_i \pm e_j \mid 1 \leq i < j \leq n \right\rbrace, \\
	\Phi^S &= \left\lbrace e_i \pm e_j \mid 1 \leq i < j \leq n \right\rbrace, \\
	\Phi^L &= \left\lbrace 2e_i \mid 1 \leq i \leq n \right\rbrace.
\end{align*}

As we have mentioned before each positive root uniquely decomposes as a linear combination
of simple roots. Let $m_i(\alpha)$ denote the coefficient at $\alpha_i$ in the decomposition of 
a root $\alpha$. Thus
\[
	\alpha = \sum_{i=1}^n m_i(\alpha) \alpha_i.
\]
It is clear, that is $\alpha, \beta, \alpha+\beta \in \Phi$ then $m_i(\alpha+\beta) = m_i(\alpha) + m_i(\beta)$.
The natural number $\sum_{i=1}^n m_i(\alpha)$ is called \textit{the height of the root $\alpha$} and is denoted by
$\height(\alpha)$. There exists the unique root in $\Phi$ of maximal height. We will denote it by $\alpha_{\max}$.
In fact, $\alpha_{\max} = 2e_1$ and $\height(\alpha_{\max}) = 2n-1$. We will occasionally use the notation
$\alpha_{ij}$ for the root $e_i-e_j$ and $\alpha_{i,-j}$ for the root $e_i+e_j$, although we 
prefer writing roots in terms of their decompositions with respect to simple ones. It's easily verified that
\begin{align*}
	\alpha_{ij} &= \alpha_i + \dots + \alpha_{j-1}, &\height(\alpha_{ij}) &= j-i \\
	\alpha_{i,-j} &= \alpha_i + \dots + \alpha_{n-1} + \alpha_n + \alpha_{n-1} + \dots + \alpha_j,
		&\height(\alpha_{i,-j}) &= 2n+1-i-j.
\end{align*}

\paragraph{Elementary root unipotents and Steinberg relations. }

Recall that $\Sp(2n,F)$ is the classical symplectic group of rank $n$ over the field $F$
and $\Up = \Up(2n,F)$ be the unipotent radical of the Borel subgroup of $\Sp(2n,F)$, i.e.
the set of all upper unitriangular matrices in $\Sp(2n,F)$. We 
introduce now the set of elements of $\Up$ known as elementary root unipotents which in fact
generate $\Up$ as a subgroup of $\Sp(2n,F)$.

Let $\alpha = \alpha_{ij}$ be a short positive root and $\xi \in F$ ($j$ being either positive or negative). 
Set $j' = \begin{cases} j, & j>0 \\ 2n+1-j, &j<0 \end{cases}$. The element
\[
	x_\alpha(\xi) = e + \xi e_{i,j'} - \varepsilon_j \xi e_{2n+1-j', 2n+1-i},
\]
where $\varepsilon_j = \operatorname{sgn}(j) = j/|j|$, is called \textit{the elementary short root unipotent}.
Given a long positive root $\alpha = \alpha_{i,-i}$ and $\xi \in F$ we will call the element
\[
	x_{\alpha}(\xi) = e + \xi e_{i,2n+1-i}
\]
\textit{the elementary long unipotent}. It's verified by a straightforward calculation that all the elementary
root unipotents that correspond to positive roots are indeed contained in $\Up$ and satisfy the following 
identities known as Steinberg relations (for all $\xi,\zeta \in F$):
\begin{description}
    \item[\mdseries\rmfamily(R2)] $x_\alpha(\xi) x_\alpha(\zeta) = x_\alpha(\xi + \zeta)$ for all $\alpha \in \Phi$;
    \item[\mdseries\rmfamily(R3)] $[x_\alpha(\xi), x_\beta(\zeta)] = e$ whenever $\alpha+\beta \notin \Phi$;
    \item[\mdseries\rmfamily(R4)] $[x_{\alpha}(\xi), x_\beta(\zeta)] = x_{\alpha+\beta}(N_{\alpha,\beta} \xi \zeta)$ 
	    whenever $\alpha, \beta, \alpha+\beta \in \Phi^S$;
    \item[\mdseries\rmfamily(R5)] $[x_\alpha(\xi), x_\beta(\zeta)] = x_{\alpha+\beta}(N_{\alpha,\beta} \xi \zeta)$ 
	    for all $\alpha, \beta \in \Phi^S$ such that $\alpha+\beta \in \Phi^L$;
    \item[\mdseries\rmfamily(R6)] $[x_\alpha(\xi), x_\beta(\zeta)] = x_{\alpha+\beta}(N_{\alpha,\beta,1} \xi \zeta) x_{\alpha+2\beta}(N_{\alpha,\beta,2} \xi \zeta^2)$ for all $\alpha \in \Phi^L, \beta \in \Phi^s$ such that $\alpha+\beta \in \Phi$.
\end{description}
The numbering of the relations oddly starts with $2$ for the sake of compatibility. \red{Elucidate!!!}
The constants $N_{*}$ are called \textit{structure constants}. In most occurrences,
it's not important what is the sign of $N_*$, but only it's module. For short roots $\alpha$ and $\beta$
such that $\alpha+\beta$ is a long root $N_{\alpha,\beta} = \pm 2$ and all other structure constants
equal $\pm 1$. 
\red{Think if we need specific signs
\begin{align*}
	N_{\alpha_{ij},\alpha_{jk}} &= - \varepsilon_k \varepsilon_j N_{\alpha_{ij},\alpha_{k,2n+1-j}} = 1 \\	
\end{align*}
}

For each $\alpha \in \Phi$ we will call the subgroup
\[
	X_\alpha = \{ x_{\alpha}(\xi) \mid \xi \in R \}
\]
\textit{the \textup{[}elementary unipotent\textup{]} root subgroup} corresponding to $\alpha$. It's a well knows fact, that
$\Ep(2n,R)$ is generated by all $X_\alpha$, $\alpha \in \bar{\Phi}$ and $\Up$ is generated 
by $X_\alpha$, $\alpha \in \Phi$ (positive roots only).
	
\paragraph{Commutators preserving maps.}

Let $G$ be a group and $\phi : G \rightarrow G$ be a bijection that preserves Lie bracket, i.e.
\[
	\phi([x,y]) = [\phi(x),\phi(y)]
\]
for every $x,y \in G$. We will call such $\phi$ \textit{a PC-map} or \textit{a commutators preserving map}
or \textit{a map preserving commutators}. It is a well known fact that PC-maps of a group form a group themselves 
under composition of maps. We will denote this group by $\PC(G)$. It is easy to see that any such map
preserves the identity element of $G$ and maps the center of $G$ bijectively onto itself. More generally,
if $X$ is a subset of $G$, then the centralizer of $X$ in $G$ is mapped bijectively onto the 
centralizer of the image of $X$ by any $\PC$-map. We will call a subgroup $H$ of $G$ stable under any
$\PC$-map a $\PC$-subgroup. Thus, the centralizer of a $\PC$-subgroup is a $\PC$-subgroup itself.

\section{Standard PC-maps of Up}
\label{sec:standard}

In this section we define some standard $\PC$-maps of $\Up$, most of which are automorphisms.
\paragraph{Inner automorphisms.} Every matrix $C \in \Up$ defines an \textit{inner automorphism}
$\Int_C : \Up \longrightarrow \Up$ as follows:
\[
	\Int_C(a) = {}^C a = C a C^{-1}
\]
for all $a \in \Up$. Clearly, inner automorphisms form a subgroup in $\PC(C)$, but we 
will see that it is not normal. \red{More on that}

\paragraph{Diagonal automorphism.} Let $d$ be a diagonal matrix in $\Sp(2n,F)$. An inner automorphism
$D = \Int_d$ of $\Sp(2n,F)$ leaves $\Up$ invariant and thus induces an automorphism of $\Up$ which 
we will call a \textit{diagonal automorphism} of $\Up$.

\paragraph{Semi-diagonal automorphism.}
Let $\epsilon \in F^*$. Define $\psi_\epsilon$ to be the map sending each $a = \begin{pmatrix}
		x & y \\
		0 & z
	\end{pmatrix} \in \Up$ to $a = \begin{pmatrix}
		x & \epsilon \cdot y \\
		0 & z
	\end{pmatrix},$ where the matrices are written with respect the partition $2n=n+n$. Clearly,
$\psi_\epsilon$ is an automorphism and if $\epsilon$ is a square in $F$ then $\psi_\epsilon$ is
a diagonal automorphism. We will call such 
automorphisms \textit{semi-diagonal}. In \cite{OuPCBSp2012, OuWangAutomBSp2008} these are called
\textit{extremal}. This, however not only is not suggestive, but also creates a collision with the 
existing terminology. Namely, in \cite{GibbsAutomUni1970} John Gibbs defines another class of automorphisms, 
which we also make use of in this paper, extremal. 
	
\paragraph{Field automorphism.} Let $\varphi : F \longrightarrow F$ be a field automorphism of $F$.
Then $\varphi$ induces an automorphism $\Sp(\varphi)$ of $\Sp(2n,R)$ by mapping a matrix $a$ to $(\varphi(a)_{ij})_{i,j=1}^n$.
Further, the restriction of $\varphi$ to $\Up$ is an automorphism of $\Up$. We will call it
a \textit{field automorphism} of $\Up$ induced by $\varphi$.

\paragraph{Extremal automorphisms.} Suppose $6F = F$. Let $u \in F$. Define an automorphism 
$\es^1_u : \Up \longrightarrow \Up$ on generators of $\Up$ as follows 
\[
	\es^1_u(x_\alpha(\xi)) = \begin{cases}
		x_{\alpha_1}(\xi) x_{\alpha_{\max}-\alpha_1}(u \xi) x_{\alpha_{\max}}(\frac{1}{2}N_{\alpha_{\max}-\alpha_1,\alpha_1} u \xi^2),
		& \alpha=\alpha_1 \\
		x_{\alpha}(\xi), & \text{ otherwise. }
	\end{cases}
\]
Further, define $\es^2_u$ as follows
\[
	\es^2_u(x_\alpha(\xi)) = \begin{cases}
		\begin{aligned}
			x_{\alpha_1}(\xi) \cdot x_{\alpha_{\max}-2\alpha_1}(u \xi) 
			\cdot x_{\alpha_{\max}-\alpha_1}(\frac{1}{2}N_{\alpha_{\max}-2\alpha_1,\alpha_1,1} u \xi^2)& \\
			\cdot x_{\alpha_{\max}}(\frac{1}{3} N_{\alpha_{\max}-2\alpha_1,\alpha_1,2} u \xi^3)&
		\end{aligned},
		& \alpha=\alpha_1 \\
		x_{\alpha}(\xi), & \text{ otherwise. }
	\end{cases}
\]

\noindent A direct calculation shows that $\es^1_u$ and $\es^2_u$ are automorphisms of $\Up$ for all $u \in F$.
Following \cite{GibbsAutomUni1970} we call $\es^1_u$ and $\es^2_u$ the \textit{extremal automorphisms}
of $\Up$ corresponding to $u$ of the first and the second type, respectively.

\paragraph{Central PC-maps} Let $\phi \in \PC(\Up)$. We will call $\phi$ a \textit{central $\PC$-map}
if $\phi(a) a^{-1} \in \Center(\Up)$ for all $a \in \Up$, i.e. $\phi$ is a multiplications by a central element map. 

\begin{lemma}
	\label{lemma:standard:central:map}
	Assume $n \geq 2$.
	Let $f : F^n \rightarrow F$ be a function such that $f(0)=0$. Define $\phi : \Up \rightarrow \Up$ as follows:
	\[
		\phi(a) = a \cdot x_{\alpha_{\max}}(f(a_{12},a_{23},\dots,a_{n,n+1})).
	\]
	Then $\phi \in \PC(\Up)$.
	\begin{proof}
		It is quite obvious that $\phi$ preserves the Lie product. Indeed, let $a,b \in \Up$. Then
		\begin{align*}
			\phi([a,b]) &= [a,b] \cdot x_{\alpha_{\max}}(0) = [a,b]\\
			&= 
			[a x_{\alpha_{\max}}(f(a_{12},a_{23},\dots,a_{n,n+1})),
			b x_{\alpha_{\max}}(f(b_{12},b_{23},\dots,b_{n,n+1}))]\\
			&=
			[\phi(a),\phi(b)].
		\end{align*}
		Next, we will  show that $\phi$ is injective. Let $a, b \in \Up$. Suppose that
		$\phi(a)=\phi(b)$. Then $a \equiv b \mod \Center(\Up)$. In particular,
		$a_{i,i+1} = b_{i,i+1}$ for $1 \leq i \leq n$. Thus,
		\begin{align*}
			a &= b \cdot x_{\alpha_{\max}}(f(b_{12},b_{23},\dots,b_{n,n+1})) \cdot 
				x_{\alpha_{\max}}(-f(a_{12},a_{23},\dots,a_{n,n+1}))\\
			&= 
			b \cdot x_{\alpha_{\max}}(f(b_{12},b_{23},\dots,b_{n,n+1})-f(a_{12},a_{23},\dots,a_{n,n+1}))\\
			&=
			b \cdot x_{\alpha_{\max}}(0) = b.
		\end{align*}
		Finally, $\phi$ is surjective. Indeed, let $b \in \Up$. Then the matrix
		\[
			a = b \cdot x_{\alpha_{\max}}(-f(b_{12},b_{23},\dots,b_{n,n+1}))
		\]
		satisfies $a_{i,i+i}=b_{i,i+1}$ for all $1 \leq i \leq n$. Thus,
		\[
			b = a \cdot x_{\alpha_{\max}}(-f(a_{12},b_{23},\dots,a_{n,n+1})) = \phi(a).
		\]
	\end{proof}
\end{lemma}
We will call the map $\phi$ constructed in the last lemma a \textit{standard central map} defined
by $f$.

\red{Talk about various opinions on this question.}

\section{Some PC-subgroups of Up.}

Recall that $\Up$ stands for the subgroup of unitriangular matrices in $\Sp(2n,R)$. 
The center $\Center(\Up)$ of $\Up$ is $X_{\alpha_{\max}} \cdot X_{\alpha_{\max}-\alpha_1}(\Ann(2))$,
where $\Ann(2) = \{\xi \in F \mid 2\xi =0 \}$ and  $X_{\alpha_{\max}-\alpha_1}(\Ann(2)) = 
\{ x_{\alpha_{\max}-\alpha_1}(\xi) \mid \xi \in \Ann(2) \}$. As in all our results we assume $2F=F$,
we have $\Ann(2)=0$ and $\Center(\Up) = X_{\alpha_{\max}}$.

Let $\Up^{(1)} = \Up$ and set for all $i \in \mathbb{N}$
\[
	\Up^{(i)} = [\Up^{(i-1)}, \Up^{(1)}].
\]
It is easily checked that $\Up^{(i)}$ is generated by all elementary root unipotents $x_{\alpha}(\xi)$
such that $\height(\alpha) \geq i$ and $\xi \in R$. In particular, $\Center(\Up)=\Up^{(2n+1)}$. 
Also,
\[
	[\Up^{(i)}, \Up^{(j)}] \subseteq \Up^{(i+j)}
\]
for all $i,j$.

Recall that al elementary root unipotents contained in $\Up$ are modulo sign symmetric with
respect to the skew-diagonal, i.e. $(x_{\alpha}(\xi))_{ij} = \pm (x_{\alpha}(\xi))_{2n+1-j,2n+1-i}$.
An arbitrary matrix $a$ in $\Up$ despite being a product of elementary root unipotents doesn't inherit
this property in general. However it holds for some entries of $a$ as stated in the following lemma.
\red{Introduce $\delta, \varepsilon$, define symplectic group in this manner}
\begin{lemma}
	\label{lemma:symmetry:of:zeros}
	Let $a \in \Up$ and $i,s \in \mathbb{N}$ be such that
	$1 \leq i \leq 2n-1$ and $1 \leq i+s \leq 2n$. Suppose that
	either $a_{i,i + k} = 0$ for all $1 \leq k \leq s-1$
	or $a'_{i+k,i+s} = 0$ for all $1 \leq k \leq s-1$.
	Then $a_{i,i + s} = - \varepsilon_i \varepsilon_{-(i + s)} a_{2n+1-(i + s),2n+1-i}$.
	
	In particular, if $a \in \Up^{(s-1)}$ then $a_{i,i + s} = \pm a_{2n+1-(i + s),2n+1-i}$.
	\begin{proof}
		First of all 
		\begin{equation}
			\label{eq:lemma:symmetry:of:zeros:1}
			0 = \delta_{i,i + s} = (a a^{-1})_{i,i + s} = 
			\sum_{j=1}^{2n} a_{ij} a'_{j,i + s}.
		\end{equation}
		Now, as $a, a^{-1} \in \Up$ it follows that $a_{ij} = 0$ for all $j<i$ and $a'_{j,i + s} = 0$ for
		all $j > i + s$. Further, by assumption of the lemma either $a_{i,i + k} = 0$ for all $1 \leq k \leq s-1$
		or $a'_{i+k,i+s} = 0$ for all $1 \leq k \leq s-1$. In any case, 
		the right-hand side of \eqref{eq:lemma:symmetry:of:zeros:1} contains only two nonzero summands and 
		we can rewrite \eqref{eq:lemma:symmetry:of:zeros:1} as 
		\[
			0 = a_{ii} a'_{i,i + s} +  a_{i, i + s} a'_{i + s, i + s} 
			= a'_{i, i + s} + a_{i,i + s} = \varepsilon_i \varepsilon_{i + s} a_{-(i + s),-i}
			+ a_{i,i + s}.
		\]
	\end{proof}
\end{lemma}

We are going to show, that $\PC$-maps are monotone in the sense, that $\Up^{(i)}$ maps bijectively onto itself
under any commutators preserving maps. The next lemma handles a particular case of this result,
namely it shows that the image of any elementary root unipotent is a product of elementary root unipotents
of (non strictly) greater height. 

\red{The assumption $2F=F$ is crucial, for otherwise not every long root unipotent is a commutator}

\begin{lemma}
	\label{lemma:transvections:in:Up:s}
	Suppose $n \geq 3$ and $2F=F$.
	Let $\phi \in \PC(\Up)$. Then $\phi(X_{\alpha}) \subseteq \Up^{(\height(\alpha))}$
	for all $\alpha \in \Phi$.
	\begin{proof}
		We will first prove by induction on height that the assertion of the lemma holds for all roots
		$\alpha \in \Phi \setminus \{ \alpha_n+\alpha_{n-1}, \alpha_n+2\alpha_{n-1} \}$.
		The base of induction, $\alpha \in \Pi$ is
		trivial. Now, suppose the assertion on the lemma holds of all roots in
		$\Phi \setminus \{ \alpha_n+\alpha_{n-1}, \alpha_n+2\alpha_{n-1} \}$ of height
		not greater than $k < 2n+1$. Let $\alpha$ be a root in 
		$\Phi \setminus \{ \alpha_n+\alpha_{n-1}, \alpha_n+2\alpha_{n-1} \}$ of height $k+1$.
		
		It can be easily seen from the construction of root
		system of type $C_n$ in Section \ref{sec:prelim} that
		if $\alpha \neq \alpha_n+\alpha_{n-1},\alpha_n+2\alpha_{n-1}$, then $\alpha$
		is a sum of short roots with the same restriction,  i.e. 
		$\alpha = \alpha'+\alpha''$ and $\alpha',\alpha'' \in \Phi^S \setminus \{ \alpha_n+\alpha_{n+1} \}$
		(note that $\alpha_n+2\alpha_{n-1}$ is a long root and thus need not be excluded here). 
		\red{Add explicit proof!}		
		It follows from the uniqueness of decomposition of roots in terms of simple roots that 
		$\height(\alpha) = \height(\alpha')+\height(\alpha'')$.	
		Note that $N_{\alpha',\alpha''} = \pm 1$ or $\pm 2$ is invertible.
		It follows that
		\[
			\phi(x_{\alpha}(\xi)) = \left[ \phi(x_{\alpha'}(\xi)),\phi\left(x_{\alpha''}\left(\frac{1}{N_{\alpha',\alpha''}}\right)\right)\right]
			\in \left[\Up^{(\height(\alpha'))}, \Up^{(\height(\alpha''))}\right] \subseteq \Up^{(\height(\alpha))}.
		\]
		
		The case of $\alpha = \alpha_n+\alpha_{n-1}$ is a bit trickier for it is not a sum
		of short roots and thus $x_{\alpha}(\xi)$ is not a commutator of elementary 
		root unipotents. Write $a = \phi(x_{\alpha}(\xi))$ as
		\begin{align*}
			\phi(x_{\alpha}(\xi)) &= b \cdot c, &b &= \prod_{i=1}^n x_{\alpha_i}(\zeta_i), &c &= \prod_{\alpha \in \Phi, \height(\alpha)\geq 2} x_{\alpha}(\zeta_\alpha).
		\end{align*}
		Suppose, $b \neq e$. Let $i$ be the minimal index such that $\zeta_i \neq 0$.
		Let $\beta = \alpha_{\max}-\alpha_1-\dots-\alpha_i$.
		$\gamma = \alpha_1+\dots+\alpha_{i-1}$ (if $i=1$ then $\gamma=0$).
		Note that as $n\geq 4$ neither $\beta$ nor $\gamma$ equals $\alpha_n+\alpha_{n-1}$ or
		$\alpha_n+2\alpha_{n-1}$. Set 
		\[
			u = [[a, x_{\beta}(1)],x_\gamma(1)],
		\] if $\gamma\neq 0$ and 
		\[	
			u = [a, x_{\beta}(1)] ,
		\]
		otherwise. A direct calculation shows that 
		\[
			u = x_{\alpha_{\max}}(\pm \zeta_i) \neq e.
		\]
		On the other hand, if $\gamma \neq 0$, then
		\begin{align*}
			\phi^{-1}(u) &= 
				[[x_{\alpha}(\xi), \phi^{-1}(x_{\beta}(1))], \phi^{-1}(x_{\gamma}(1))] \\
				&\subseteq 
				[[\Up^{(2)}, \Up^{(2n+1-i)}], \Up^{(i-1)}] \subseteq \Up^{(2n+2)} = \{e\};
		\end{align*}
		and if $\gamma = 0$, then
		\begin{align*}
			\phi^{-1}(u) = 		
				[x_{\alpha}(\xi), \phi^{-1}(x_{\beta}(1))] \subseteq 
				[\Up^{(2)}, \Up^{(2n)}] \subseteq \Up^{(2n+2)} = \{e\}.
		\end{align*}
		In both cases we get a contradiction with the fact that 
		any $\PC$-map preserves the neutral element. Thus, $b=0$ and 
		$\phi(x_{\alpha}(\xi)) \in \Up^{(2)}$. Finally, the case if $\alpha = \alpha_n+2\alpha_{n-1}$
		then
		\[
			\phi(x_{\alpha}(\xi)) = \left[ \phi(x_{\alpha_{n-1}}(\xi)),\phi\left(x_{\alpha_n+\alpha_{n-1}}\left(\frac{1}{N_{\alpha_{n-1},\alpha_n+\alpha_{n-1}}}\right)\right)\right]
			\in \left[\Up, \Up^{(2)}\right] \subseteq \Up^{(3)}.
		\]		
	\end{proof}		
\end{lemma}

\red{Intro for the next three lemmas}

\red{Make sure $\U_1$ is defined}
\red{In these lemmas the base 1,..,-1 would be more convenient}
\begin{lemma}
	\label{lemma:extract:from:U1}
	Suppose $n \geq 1$. Let $a \in \U_1$. Then for any $2 \leq j \leq 2n-1$
	there exist a root $\beta$ of height $2n-j$ such that
	\[
		x_{\alpha_{\max}}(\pm a_{1j}) = [x_{\beta}(1), a].
	\]
	\begin{proof}
		Note that
		\[
			a = x_{\alpha_{\max}}(*) \prod_{\alpha \in \Phi, m_1(\alpha)=1} x_{\alpha}(a_{1j}).
		\]
		Set $\beta = \alpha_{\max}-\alpha_{1j}$. Note that the only root $\alpha \in \Phi$ 
		such that $m_1(\alpha)\geq 1$ and $\alpha + \beta$ is a root is $\alpha_{1j}$. Thus
		\[
			[x_{\beta}(1),a] = x_{\alpha_{\max}}(\pm 2 a_{1j}).
		\]
		Finally, note that $\height(\beta) = 2n-j$.
	\end{proof}
\end{lemma}

\begin{lemma}
	\label{lemma:extract:middle}
	Suppose $n \geq 2$. Fix a pair of indices $2 \leq i<j \leq 2n-1$.
	Then there exist roots $\beta,\gamma \in \Phi$ such that
	\[
		x_{\alpha_{\max}}(\pm 2 a_{ij}) = [x_\beta(1),[x_\gamma(1),a]]
	\]
	for any matrix $a$ in $\Up$. Further, $\height(\beta)+\height(\gamma) = 2n-1+i-j$.
	\begin{proof}
		Set $j' = -j$ if $j\leq n$ and $j' = 2n+1-j$, otherwise.
		By assumption $j \neq 1, 2n$, thus 
		$\gamma = \alpha_{1,j'}$ is a short root. Let
		\begin{align*}
			b &= [x_\gamma(1),a]
			&
			\tilde{b} &= a x_\gamma(-1) a^{-1}.
		\end{align*}				
		By a direct calculation we see that
		\[
			\tilde{b} = e + \pm a_{*1} a'_{2n+1-j,*} \pm a_{*j} a'_{2n,*} = e \pm e_{*1} a'_{2n+1-j,*} \pm a_{*j} e_{2n,*}.
		\]
		Note that $\tilde{b}$ and thus also $b$ are contained in $\U_1$.
		Also, $b_{1,2n+1-i} = \tilde{b}_{1,2n+1-i} = \pm a'_{2n+1-j,2n+1-i} = \pm a_{ij}$.
		By Lemma \ref{lemma:extract:from:U1} there exists a root $\beta$ of height $i-1$
		such that $x_{\alpha_{\max}}(\pm 2 a_{ij}) = [x_\beta(1),b]$. It's only left
		to notice that $\height(\gamma) = j-1$. Thus, indeed $\height(\beta)+\height(\gamma) = 2n-1+i-j$.
	\end{proof}
\end{lemma}

\begin{lemma}
	\label{lemma:PC:preserves:Up:s}
	Assume $n \geq 3$ and $2F=F$.
	Let $\phi \in \PC(\Up)$. Then $\phi(\Up^{(s)}) = \Up^{(s)}$
	for all $s \leq 2n-1$.
	\begin{proof}
		By invertibility of $\phi$ it is enough to show that  $\Up^{(s)}) \subseteq \Up^{(s)}$ for all
		$1 \leq s \leq 2n-1$. This obviously holds for $s=1$ as well as $s = 2n-1$.
		Suppose the lemma holds for all $s' \leq s$, where $s<2n-2$. We will prove that it holds for $s+1$.
		Pick an arbitrary matrix $x$ in $\Up^{(s+1)}$. Suppose that $a = \phi(x) \notin U^{(s+1)}$. By the assumption
		of induction $a$ lies in $\Up^{(s)} \setminus \Up^{(s+1)}$. By Lemma \ref{lemma:symmetry:of:zeros}
		this yields that there exists such an index $1 \leq i \leq n$ such that $a_{i,i + s} \neq 0$. 
		Note that $1 \leq i < i+s \leq 2n$. If $i>1$ then by Lemma \ref{lemma:extract:middle}
		there exist roots $\beta, \gamma \in \Phi$ such that
		\[
			x_{\alpha_{\max}}(\pm 2 a_{i,i+s}) = [x_\beta(1),[x_\gamma(1),a]]
		\]
		and $\height(\beta)+\height(\gamma) = 2n-1-s$.
		By Lemma \ref{lemma:transvections:in:Up:s} we have 
		\[
			\phi^{-1} (x_{\alpha_{\max}}(\pm 2 a_{i,i+s})) =
			[\phi^{-1}(x_\beta(1)),[\phi^{-1}(x_\gamma(1)),x]] \in
			[\Up^{(\height(\beta))},[\Up^{(\height(\beta))},\Up^{(s+1)}]] \le
			\Up^{(2n)} = e.
		\]
		Thus, $x_{\alpha_{\max}}(\pm a_{i,i+s}) = e$ and consequently $a_{i,i+s}=0$, which contradicts with the
		assumption that $a_{i,i+s}\neq 0$.
		
		Now suppose $a_{i,i+s} = 0$ for all $1 < i \leq n$ and $a_{1,1+s} \neq 0$. Then
		\[
			a = x_{\alpha_{1j}}(a_{1,1+s}) \cdot \prod_{\alpha \in \Phi, \height(\alpha)\geq s+1} x_{\alpha}(a_\alpha),
		\]	
		where $j = 1+s$, if $1+s \leq n$ and $j = 2n-s$, otherwise. It's clear that 
		\[
			[x_{\alpha_{\max}-\alpha_{1,j}}(1),a] = x_{\alpha_{\max}}(\pm 2 a_{1j}),
		\]
		as $\height(\alpha_{\max}-\alpha_{1,j})=2n-1-s$ and thus the sum of the root 
		$\alpha_{\max}-\alpha_{1,j}$ and a positive root of height not less than $s+1$ is never a root.
		Thus, again by Lemma \ref{lemma:transvections:in:Up:s},
		\[
			\phi^{-1}(x_{\alpha_{\max}}(\pm 2 a_{1j})) = 
				[\phi^{-1}(x_{\alpha_{\max}-\alpha_{1,j}}(1)),x] \in \Up^{(2n)} = \{e\}.
		\]
		it follows that $a_{1j} = 0$.
	\end{proof}
\end{lemma}

The definition of the subgroups $P^i_k$ mimic the definition of the subgroups $P_k$ of \red{REF}.
There, these subgroups appear purely as a technical gadget and the proof that they are $\PC$-subgroups
is carried out in the same way as the proof that the members of upper/lower central series are
$\PC$-subgroups. However, we will show that the subgroup $P^i_k$ and the subgroups $P_k$ as the particular
case with $i=1$, arise naturally as the intersections of the centralizers of some known $\PC$-subgroups
with some other known $\PC$-subgroups.

\begin{lemma}
	\label{lemma:PC:preserves:P:i:k}
	Let $\phi \in \PC(\Up)$. Then $\phi$ preserves the sets
	\[
		P^i_k = \Up^{(i+1)} \cdot \prod_{\alpha \in \Phi, \height(\alpha) = i, m_1(\alpha) = \dots = m_{k-1}(\alpha) = 0} x_\alpha(\xi_\alpha)
	\] for all $1 \leq i \leq 2n-1$ and $1 \leq k \leq \lceil\frac{2n-i}{2}\rceil.$ \red{Check the upper bound}.
	\begin{proof}
		\red{Explicate:}
		It is easily shown that
		\[
			P^i_k = C_{\Up}(\Up^{(2n-k)}) \cap \Up^{(i)}
		\] and the right-hand side is preserved by $\phi$ by
		Lemma \ref{lemma:PC:preserves:Up:s}. Thus $P_k^i$ is preserved by $\phi$ as well.
	\end{proof}
\end{lemma}

The last lemma is extremely useful via the following observation. Let $\phi$ be a $\PC$-map and $x_\alpha(\xi)$ 
an elementary root unipotent. So far we don't know already quite some data about the image of $x_\alpha(\xi)$ under $\phi$.
First, by Lemma \ref{lemma:transvections:in:Up:s} we have $\phi(x_\alpha(\xi) \in \Up^{(\height(\alpha))}$.
Then
\[
	\phi(x_\alpha(\xi)) = \prod_{\beta \in \Phi, \height(\beta)=\height(\alpha)} x_{\beta}(\zeta_\beta) \cdot x, 
	\qquad x \in \Up^{(\height(\alpha)+1)}.
\]
Let $k$ be the least number such that $m_k(\alpha)>0$. Then $x_{\alpha}(\xi) \in P^{(\height(\alpha))_k} \setminus
P^{(\height(\alpha))_{k+1}}$ \red{extreme case $k=n$}.
By Lemma \ref{lemma:PC:preserves:P:i:k} it follows that $\zeta_\beta = 0$ whenever $m_i(\beta)>0$ for some
$i<k$. In particular, if $\alpha = \alpha_{ij}$ and $j' = \begin{cases} j, &j>0 \\ 2n+1+j, &j<0 \end{cases}$
(and thus $(x_\alpha(\xi))_{i,j'} = \xi$) we have $(\phi(x_\alpha(\xi)))_{i,j'} \neq 0$.

\[
	\red{U_1^l = \Up^{(l)} \cdot X_{\alpha_1} \cdot X_{\alpha_1 + \alpha_2} \cdot \dots \cdot X_{\alpha_1 + \dots +\alpha_{l-1}}}
\]

\section{Computing centralizers}

In this section we will develop the technical framework for computing centralizers of products of root subgroups
of $\Up$ and compute certain centralizers for the later use in Section \ref{sec:up:to:ai}. The following lemma is a
simple but powerful tool for reducing the computation of centralizers to combinatorics of roots.

\red{add a mention about the case $x_{\alpha_{\max}}(*))$ to the proof}
\begin{lemma}
	\label{lemma:centralizer:X}
	Suppose $2F=F$. Let $\alpha \in \Phi$. Then 
	\[
		C_G(X_\alpha) = C_G(x_\alpha(\xi)) = C_G(x_\alpha(\xi)x_{\alpha_{\max}}(*)) = \prod_{\beta \in \Phi, \alpha+\beta \notin \Phi} X_\beta
	\]
	for any $\xi \neq 0$.
	\begin{proof}
		Set $\alpha^\perp = \{\beta \in \Phi \mid \alpha+\beta \notin \Phi\}$.
		Consider an element $a = x_\alpha(\xi) \in X_\alpha\setminus\{e\}$. It follows form the Steinberg
		relations that
		\[
			\prod_{\beta \alpha^\perp} X_\beta \subseteq C_G(x_\alpha(\xi)).
		\]
		Now let $b \notin \prod_{\beta \in \alpha^\perp} X_\beta$. We will show
		that $b$ doesn't centralize $a$. Indeed, let
		\begin{equation}
			\label{lemma:centralizer:X:1}
			\prod_{\beta \in \alpha^\perp} x_\beta(\zeta_\beta).
		\end{equation}
		Choose a root $\beta' \in \alpha^\perp$ such that $\beta'+\alpha \notin \Phi$ and $\zeta_{\beta'} \neq 0$.
		Moreover, such $\beta'$ can be chosen minimal \red{(with respect to the root the set of positive roots $\Pi$)}
		among roots in $\beta'$. Then $\alpha+\beta$ is not an integral linear combination of 
		$\{\alpha\} \cup (\alpha^\perp \setminus \{\beta'\})$. Without loss of generality, we may assume that
		$\beta'$ goes first in the product \eqref{lemma:centralizer:X:1}. It follows that
		\[
			[a,b] = [x_\alpha(\xi), x_{\beta'}(\zeta_{\beta'}) \cdot 
			\prod_{\beta \in \alpha^\perp, \beta \neq \beta'} x_\beta(\zeta_\beta)]
			= [x_{\alpha+\beta'}(N_{\alpha,\beta',1} \xi \zeta_{\beta'}) \cdot 
			\prod_{\gamma \in S} x_\gamma(*)],
		\]
		where $S \subseteq \langle \{\alpha\} \cup (\alpha^\perp \setminus \{\beta'\}) \rangle \cap \Phi$.
		Note that $\alpha+\beta'$ is minimal in the set $\{\alpha+\beta'\}\cup S$ and thus $[a,b] = e$ if and only if
		$N_{\alpha,\beta',1} \xi \zeta_{\beta'} = 0$. Thus $\zeta_{\beta'} = 0$.
%
%
%
	\end{proof}
\end{lemma}

\red{Check if any dependent on the rank}
\begin{lemma}
	\label{lemma:simple:roots:as:C}
	Suppose $2F = F$. Then 
	\begin{enumerate}
		\item for all $1 \leq i \leq n$ the centralizer of
			\[
				\prod_{\alpha \in \Phi, m_1(\alpha)\geq 1, \alpha \neq \alpha_1+\dots+\alpha_{i-1}, 
					\alpha \neq \alpha_{\max}-\alpha_1-\dots-\alpha_i} X_{\alpha},
			\]
			equals
			\[
				X_{\alpha_i} \cdot X_{\alpha_1+\dots+\alpha_i} \cdot 
				X_{\alpha_{\max}-\alpha_1-\dots-\alpha_{i-1}} \cdot X_{\alpha_{\max}};
			\]
		\item the centralizer of 
		\[
			X_{\alpha_1} X_{\alpha_2} \dots X_{\alpha_{n-3}} \cdot
			X_{\alpha_{\max}-\alpha_1-\dots-\alpha_{n-1}} \cdots
			X_{\alpha_{\max}-\alpha_1} \cdot			
			X_{\alpha_{\max}}
		\]
		equals
		\[
			X_{\alpha_n} \cdot X_{\alpha_{n-1}+\alpha_n} \cdot X_{2\alpha_{n-1}+\alpha_n} \cdot 
			X_{\alpha_1+\dots+\alpha_n} \cdot X_{\alpha_1+\dots+\alpha_n +\alpha_{n-1}} X_{\alpha_{\max}};
		\]
		\item the centralizer of 
		\[
			\prod_{\alpha \in \Phi, m_1(\alpha)\geq 1, 
					\alpha \neq \alpha_{\max}-\alpha_1, 
					\alpha \neq \alpha_{\max}-\alpha_1-\alpha_2} X_{\alpha},
		\]
		equals
		\[
				X_{\alpha_1} \cdot 	X_{\alpha_1+\alpha_2} \cdot X_{\alpha_{\max}};
		\]
		\item the centralizer of 
		\[
			\prod_{\alpha \in \Phi, m_1(\alpha)\geq 1, 
					\alpha \neq \alpha_1, 
					\alpha \neq \alpha_1+\alpha_2,
					\alpha \neq \alpha_1+\dots+\alpha_{n-1}} X_{\alpha},
		\]
		equals
		\[
			\begin{aligned}
				X_{\alpha_n} &\cdot
				X_{\alpha_2+\dots+\alpha_n}	\cdot
				X_{\alpha_2+\dots+\alpha_n+\dots+\alpha_3} \cdot
				X_{\alpha_3+\dots+\alpha_n}	\cdot
				X_{\alpha_2+\dots+\alpha_n+\dots+\alpha_2} \\
				&\cdot
				X_{\alpha_3+\dots+\alpha_n+\dots+\alpha_3} \cdot
				X_{\alpha_{\max}-\alpha_1} \cdot
				X_{\alpha_{\max}-\alpha_1-\alpha_2} \cdot
				X_{\alpha_1+\dots+\alpha_n} \cdot
				X_{\alpha_{\max}}.
			\end{aligned}
		\]
	\end{enumerate}
	\begin{proof}
		Recall that the centralizer of a product of subgroups is the intersection of the centralizers of factors.
		Let $1 \leq i < n$. By Lemma \ref{lemma:centralizer:X} it is enough to check that
		\begin{align*}
			S &= \bigcap_{\alpha \in \Phi, m_1(\alpha)\geq 1, \alpha \neq \alpha_1+\dots+\alpha_{i-1}, 
					\alpha \neq \alpha_{\max}-\alpha_1-\dots-\alpha_i} \alpha^\perp \\
			&= \{ \alpha_i, \alpha_1+\dots+\alpha_i, \alpha_{\max}-\alpha_1-\dots-\alpha_{i-1}, \alpha_{\max}\},
		\end{align*}
		which is checked straightforwardly using the explicit description of roots of $C_n$ given in Section
		\ref{sec:prelim}. The rest of the lemma is handled in the same way.
	\end{proof}
\end{lemma}

The last lemma, despite being nicely formulated is not very useful in this form for our goals. The reason
being that we rarely can show that the subgroups, whose centralizes are computed in Lemma \ref{lemma:simple:roots:as:C},
are bijectively mapped onto themselves by a $\PC$-map in consideration. Usually, we can trace only
the images of some elementary root unipotents. Luckily, this turns out to be sufficient, as shown in the following corollary.
\begin{corollary}
	\label{cor:centralizers}
	Assume $2F = F$. Let $\phi \in \PC(\Up)$ and $S \subseteq \Phi$. Suppose that
	\[
		\phi(x_\alpha(1)) \in X_\alpha \cdot X_{\alpha_{\max}}, \quad \text{ for all } \alpha \in S.
	\]
	Then the centralizer $C_{\Up}(\prod_{\alpha \in S} X_{\alpha})$
	is preserved by $\phi$.
	\begin{proof}
		First, not that by Lemma \ref{lemma:centralizer:X} holds
		\begin{equation}
			\label{eq:cor:centralizers:1}
			C_{\Up}(\prod_{\alpha \in S} X_{\alpha}) =
			\bigcap_{\alpha \in S} C_{\Up}(X_\alpha) = 
			\bigcap_{\alpha \in S} C_{\Up}(x_\alpha(1)).
		\end{equation}
		Next, $\phi(x_\alpha(1)) = x_\alpha(\zeta_\alpha) x_{\alpha_{\max}}(*)$.
		By \red{the corollary of } Lemma \ref{lemma:PC:preserves:P:i:k}, $\zeta_\alpha \neq 0$,
		thus by Lemma \ref{lemma:centralizer:X} we have
		\begin{equation*}
			C_{\Up}(x_\alpha(1)) = C_{\Up}(\phi(x_\alpha(1))).
		\end{equation*}
		Recall, that $\PC$-maps map centralizers of elements of $\Up$ to the centralizers of their images,
		thus
		\begin{equation}
			\label{eq:cor:centralizers:2}
			\phi(C_{\Up}(x_\alpha(1))) = C_{\Up}(x_\alpha(\zeta_\alpha) x_{\alpha_{\max}}(*)) = 
			C_{\Up}(x_\alpha(1))
		\end{equation}
		for all $\alpha \in S$. Combining \eqref{eq:cor:centralizers:1}, \eqref{eq:cor:centralizers:2} 
		and the fact that $\phi$ is bijective we get
		\begin{align*}
			\phi(C_{\Up}(\prod_{\alpha \in S} X_{\alpha})) &=
			\phi(\bigcap_{\alpha \in S} C_{\Up}(x_\alpha(1))) =
			\bigcap_{\alpha \in S} \phi(C_{\Up}(x_\alpha(1)))\\
			&=
			\bigcap_{\alpha \in S} C_{\Up}(x_\alpha(1)) =
			C_{\Up}(\prod_{\alpha \in S} X_{\alpha}).
		\end{align*}
	\end{proof}
\end{corollary}

\section{Up to an almost identity map}
\label{sec:up:to:ai}

\begin{lemma}
	\label{lemma:T12}
	Assume $\rk(\Phi)\geq 3$, i.e. $n \geq 3$.
	Let $\phi \in \PC(\Up)$. Then $\phi(X_{\alpha_1}) \le \U_1 \cdot X_{\alpha_{\max} - 2\alpha_1} \cdot X_{\alpha_{\max}-2\alpha_1-\alpha_2} \cdot X_{\alpha_{\max}-2\alpha_1-2\alpha_2}$.
	\begin{proof}
		Let 
		\[
			H^{(i)}_k = \prod_{\alpha \in \Phi, \height(\alpha)=i, m_1(\alpha) = \dots = m_{k-1}(\alpha)=0} 
			X_{\alpha} \cdot U_1 \cdot \Up^{(i+1)}.
		\]
		We will show by induction on two parameters, $i$ and $k$ that $\phi(x_{\alpha_1}(\xi)) \in H^{(i)}_k$
		for all $i \leq 2n-4$.
	
		First, we will show that if $i \leq 2n-5$ then $\phi(x_{\alpha_1}(\xi)) \in H^{(i)}_k$
		yields $\phi(x_{\alpha_1}(\xi)) \in H^{(i)}_{k+1}$. If $H^{(i)}_k = H^{(i)}_{k+1}$ there is nothing 
		to prove. In particular, we may assume straight ahead that $k \geq 2$.
		Assume that $H^{(i)}_k \neq H^{(i)}_{k+1}$. Then there exists a root $\beta \in \Phi$ such that
		$\height(\beta)=i$, $m_1(\beta)=\dots=m_{k-1}(\beta)=0$ and $m_k(\beta)>0$. In fact,
		in the root system of type $C_n$ such root is unique. Write $\phi(x_{\alpha_1}(\xi))$ as follows:
		\[
			\phi(x_{\alpha_1}(\xi)) = x_\beta(\zeta) \cdot
				\prod_{\alpha \in \Phi, \height(\alpha)=i, m_1(\alpha) = \dots = m_{k}(\alpha)=0} 
				X_{\alpha} \cdot \U_1 \cdot \Up^{(i+1)}.
		\]
		Suppose there exists a simple root $\alpha_j \in \Pi \setminus \{\alpha_1\}$ such that $\alpha+\beta \in \Phi$,
		but $\alpha+\alpha_1 \notin \Phi$. Consider the commutator $[\phi(x_{\alpha_j}(1)),\phi(x_{\alpha_1}(\xi))]$. 
		We aim to show that this commutator can equal zero only if $\zeta=0$.
		By assumption, $\phi(x_{\alpha_1}(\xi)) \in H^{(i)}_k$.
		Further, by \red{the corollary of} Lemma \ref{lemma:PC:preserves:P:i:k}
		$\phi(x_{\alpha_{j}}(1))$ can be written $x_{\alpha_{j}}(\eta) \cdot y, y \in P^1_{j+1}$
		with $\eta \neq 0$. Then it is clear that
		\begin{equation}
		\label{eq:lemma:T12:bigcomm}
			\begin{aligned}
				[\phi(x_{\alpha_{j}}(1)), \phi(x_{\alpha_1}(\xi))] 
				&\in 
				[x_{\alpha_{j}}(\eta) \cdot y, x_\beta(\zeta) \cdot
				\prod_{\alpha \in \Phi, \height(\alpha)=i, m_1(\alpha) = \dots = m_{k}(\alpha)=0} 
				X_{\alpha} \cdot \U_1 \cdot \Up^{(i+1)}] \\
				&\le [x_{\alpha_{j}}(\eta) \cdot y, x_\beta(\zeta)] \\
				&\quad\cdot
				{}^{x_\beta(\zeta)} [
				x_{\alpha_{j}}(\eta) \cdot y, \prod_{\alpha \in \Phi, \height(\alpha)=i, 
					m_1(\alpha) = \dots = m_{k}(\alpha)=0} 
				X_{\alpha}] \\
				&\quad \cdot
				{}^{x_\beta(\zeta) \cdot
				\prod_{\alpha \in \Phi, \height(\alpha)=i, m_1(\alpha) = \dots = m_{k}(\alpha)=0} 
				X_{\alpha}}
				[x_{\alpha_{j}}(\eta) \cdot y, \U_1]  \\
					&\quad\cdot
				{}^{x_\beta(\zeta) \cdot
				\prod_{\alpha \in \Phi, \height(\alpha)=i, m_1(\alpha) = \dots = m_{k}(\alpha)=0} X_{\alpha} \cdot \U_1}  
				[x_{\alpha_{j}}(\eta) \cdot y, \Up^{(i+1)}].
			\end{aligned}
		\end{equation}
		We will compute the factors of the right-hand side individually. First,
		\begin{equation}
			\label{eq:lemma:T12:bigcomm:1}
			\begin{aligned}
				[x_{\alpha_{j}}(\eta) \cdot y, x_\beta(\zeta)] 
				&= 
				{}^{x_{\alpha_{j}}(\eta)} [y, x_\beta(\zeta)] \cdot
				[x_{\alpha_{j}}(\eta), x_\beta(\zeta)] \\
				&= x_{\alpha_j+\beta}(N_{\alpha_j,\beta,1} \eta \zeta) \cdot
				\prod_{\gamma \in \Phi, \height(\gamma) \geq i+1, \gamma \neq \alpha_j+\beta} x_\gamma(*),
			\end{aligned}
		\end{equation}
		where on the right-hand side none of the root $\gamma$ equal $\alpha_j+\beta$ because all root unipotents
		comprising $y$ correspond either to simple roots other than $\alpha_j$ or roots of height at least 2. Next,
		consider
		\begin{equation}
			\label{eq:lemma:T12:bigcomm:2}
			\begin{aligned}
			{}^{x_\beta(\zeta)} [x_{\alpha_{j}}(\eta) \cdot y,& \prod_{\alpha \in \Phi, \height(\alpha)=i, 
					m_1(\alpha) = \dots = m_{k}(\alpha)=0} X_{\alpha}] \\
			&= 
			{}^{x_\beta(\zeta) \cdot x_{\alpha_{j}}(\eta)} [y, \prod_{\alpha \in \Phi, \height(\alpha)=i, 
					m_1(\alpha) = \dots = m_{k}(\alpha)=0} X_{\alpha}] \\
			&\quad\cdot
			{}^{x_\beta(\zeta)} [x_{\alpha_{j}}(\eta), \prod_{\alpha \in \Phi, \height(\alpha)=i, 
					m_1(\alpha) = \dots = m_{k}(\alpha)=0} X_{\alpha}] \\
			&\le
			{}^{x_\beta(\zeta) \cdot x_{\alpha_{j}}(\eta)} [y, \prod_{\alpha \in \Phi, \height(\alpha)=i, 
					m_1(\alpha) = \dots = m_{k}(\alpha)=0} X_{\alpha}] \\
			&\quad\cdot \prod_{\gamma \in \Phi, \height(\gamma)\geq i+1, \gamma\neq\alpha_j+\beta} x_\gamma(*).
			\end{aligned}
		\end{equation}
		Now consider the commutator
		$[x_{\alpha_{j}}(\eta), \prod_{\alpha \in \Phi, \height(\alpha)=i, m_1(\alpha) = \dots = m_{k}(\alpha)=0} X_{\alpha}]$.
		Write $y$ as $\prod_{\gamma \in \Phi, m_1(\gamma)=\dots=m_j(\gamma)=0} x_{\gamma}(\theta_\gamma)$. We are going to show
		Suppose, $j \geq k$. $m_j(\gamma+\alpha) = 0$ for all $\alpha \in \Phi, \height(\alpha)=i, 
		m_1(\alpha) = \dots = m_{k}(\alpha)=0$ and all 
		$\gamma \in \Phi: m_1(\gamma)=\dots=m_j(\gamma)=0$. Thus, $\gamma+\alpha \neq \alpha_j+\beta$
		for the same $\alpha$ and $\gamma$. Thus 
		\begin{equation}
			\label{eq:lemma:T12:bigcomm:2.5}
			[x_{\alpha_{j}}(\eta), \prod_{\alpha \in \Phi, \height(\alpha)=i, m_1(\alpha) = \dots = m_{k}(\alpha)=0} X_{\alpha}] \le \prod_{\gamma \in \Phi, \height(\gamma)=i+1, \gamma \neq \alpha_j+\beta} x_\gamma(*).
		\end{equation}
		Then \eqref{eq:lemma:T12:bigcomm:2} and \eqref{eq:lemma:T12:bigcomm:2.5} together yield that
		\begin{equation}
			\label{eq:lemma:T12:bigcomm:2'}
			\begin{aligned}
			{}^{x_\beta(\zeta)} [x_{\alpha_{j}}(\eta) \cdot y,& \prod_{\alpha \in \Phi, \height(\alpha)=i, 
					m_1(\alpha) = \dots = m_{k}(\alpha)=0} X_{\alpha}] \le
			\prod_{\gamma \in \Phi, \height(\gamma)\geq i+1, \gamma\neq\alpha_j+\beta} x_\gamma(*).
			\end{aligned}
		\end{equation}
		Finally,
		\begin{equation}
			\label{eq:lemma:T12:bigcomm:3}
			\begin{aligned}
				&{}^{x_\beta(\zeta) \cdot
				\prod_{\alpha \in \Phi, \height(\alpha)=i, m_1(\alpha) = \dots = m_{k}(\alpha)=0} 
				X_{\alpha}}
				[x_{\alpha_{j}}(\eta) \cdot y, \U_1]  \\
				&\quad\cdot
				{}^{x_\beta(\zeta) \cdot
				\prod_{\alpha \in \Phi, \height(\alpha)=i, m_1(\alpha) = \dots = m_{k}(\alpha)=0} X_{\alpha} \cdot \U_1}  
				[x_{\alpha_{j}}(\eta) \cdot y, \Up^{(i+1)}] \\
				&\le
				\U_1 \cdot \Up^{(i+2)}.
			\end{aligned}
		\end{equation}
		Collecting \eqref{eq:lemma:T12:bigcomm:1}, \eqref{eq:lemma:T12:bigcomm:2'} and \eqref{eq:lemma:T12:bigcomm:3}
		we deuce that
		\begin{equation}
			\label{eq:lemma:T12:bigcomm:result}
			[\phi(x_{\alpha_{j}}(1)), \phi(x_{\alpha_1}(\xi))] =
				x_{\alpha_j+\beta}(N_{\alpha_j,\beta,1}) \cdot
				\prod_{\gamma \in \Phi, \height(\gamma) \geq i+1, \gamma \neq \alpha_j+\beta} x_\gamma(*) \cdot \U_1.
		\end{equation}
		Note that, no two roots on the right-hand side of \eqref{eq:lemma:T12:bigcomm:result}
		give $\beta+\alpha_j$ as their sum, it follows that the matrix entry
		of $[\phi(x_{\alpha_{j}}(1)), \phi(x_{\alpha_1}(\xi))]$ at the position
		corresponding to the root $\beta+\alpha_j$ equals $N_{\alpha_j,\beta,1} \eta \zeta$.
		On the other hand, by assumption $\alpha_j$ commutes with $\alpha_1$, thus
		\[
			e = \phi([x_{\alpha_1}(\xi), x_{\alpha_j}(1)]) = [\phi(x_{\alpha_{j}}(1)), \phi(x_{\alpha_1}(\xi))].
		\]
		This yields that $\zeta = 0$. 
		\red{Make a convention that $N_{\beta,\alpha_{k-1},1} = N_{\beta,\alpha_{k-1}}$ for short $\beta$.}
		\red{Make this observation about matrix entries corresponding to minimal roots.}

		Now we will choose such $\alpha_j$ as above for each $i$ and $k$. If $k>3$ this is very easy:
		set $\alpha_j=\alpha_{k-1}$. By assumption, $k \geq 2$, so the remaining cases are $k=2,k=3$.
		Let $\beta = \alpha_{pq}$. If $j>0$ we can set $\alpha_j = \alpha_q$. If $q<0$, then $p<q$
		and thus $p=k \in \{2,3\}$. Further, $2n+1-p-|q| = \height(\alpha) = i$.
		Thus $|q| = 2n+1-k-i$. While $i+k < 2n-2$ we have $|q| > 3 $ and we can set $\alpha_j = \alpha_q$.
		Note that the condition $i+k < 2n-2$ is fulfilled whenever $i < 2n-5$ as well as
		for $i = 2n-5$ and $k=2$. Thus, we have by induction that $\phi(x_{\alpha_1}(\xi)) \in H^{(2n-5)}_3$.
		In other words,
		\[
			a = \phi(x_{\alpha_1}(\xi)) \in \U_1 \cdot X_{\alpha_{\max} - 2\alpha_1} 
				\cdot X_{\alpha_{\max}-2\alpha_1-\alpha_2}
				\cdot X_{\alpha_{\max}-2\alpha_1-2\alpha_2}.
		\]
	\end{proof}
\end{lemma}

Set $\U_1^{(1)} = U_1$. For $k\geq 2$ set Let $\U_1^{(k)} = [\U^{(k-1)}_1, \Up]$. 
Note that $\U_1^{(k)} = \U_1 \cap \Up^{(k)}$. Each $\U_1^{(i)}$ is a normal subgroup of $\Up$.
\red{Should be defined somewhere in the preliminaries.}

The next lemma is a technical result based on a simple idea: if some property doesn't hold for 
all know examples, there should be a reason for it not to hold for all objects in consideration.
On the technical side, this is probably the most detailed computation in the paper. It was
originally performed using a computer program. However it is possible to verify it by hand.

\begin{lemma}
	\label{lemma:laundry}
	Suppose $6F=F$ and $n \geq 3$. Let $\phi \in \PC(\Up)$. Suppose that
	\begin{equation}
		\label{eq:lemma:laundry:assumption}
		\begin{aligned}
			\phi(x_{\alpha_1}(1)) &= x_{\alpha_1}(a_1) \cdot x_{\alpha_{\max}}(a_m) \cdot x_{\alpha_{\max}-\alpha_1}(a_{m1})
				\cdot x_{\alpha_{\max}-2\alpha_1}(a_{m11}) \\
				&\quad\cdot
				x_{\alpha_{\max}-2\alpha_1-\alpha_2}(a_{m112}) 
				\cdot x_{\alpha_{\max}-2\alpha_1-2\alpha_2}(a_{m1122})\\
			\phi(x_{\alpha_1+\alpha_2}(1)) &\in x_{\alpha_1+\alpha_2}(c_{12}) 
				\cdot x_{\alpha_{\max}}(c_{m}) \cdot x_{\alpha_{\max}-\alpha_1}(c_{m1})
				\cdot x_{\alpha_{\max}-\alpha_1-\alpha_2}(c_{m12}) \\
				&\quad\cdot 
				x_{\alpha_{\max}-2\alpha_1-\alpha_2}(c_{m112}) \cdot
				x_{\alpha_{\max}-2\alpha_1}(c_{m11}) \cdot x_{\alpha_{\max}-2\alpha_1-2\alpha_2}(c_{m1122})\\
		\end{aligned}
	\end{equation}
	\red{Probably, the last one is not there... $c_{m1122}$ should be removed; Check the usage!}
	Then 
		\begin{equation}
		\label{eq:lemma:laundry:conclusion}
		\begin{aligned}
			\phi(x_{\alpha_1}(1)) &\in X_{\alpha_1} \cdot X_{\alpha_{\max}} \cdot X_{\alpha_{\max}-\alpha_1} 
			\cdot X_{\alpha_{\max}-2\alpha_1} \\
			\phi(x_{\alpha_1+\alpha_2}(1)) &\in X_{\alpha_1+\alpha_2} 
			\cdot X_{\alpha_{\max}} \cdot X_{\alpha_{\max}-\alpha_1} 
			\cdot X_{\alpha_{\max}-\alpha_1-\alpha_2}.
		\end{aligned}
	\end{equation}
	\begin{proof}
		The proof is essentially very simple. Let $a = \phi(x_{\alpha_1}(1))$, $b = \phi(x_{\alpha_2}(1))$ and
		$c=\phi(x_{\alpha_1+\alpha_2}(1))$. Then, as $\phi$ preserves commutators, it follows that
		\begin{align*}
			\label{eq:lemma:laundry:comms}
			[a,b]&=c, &[a,c]&=e, &[c,b]&=e.
		\end{align*}
		
		Computing the commutators explicitly we get a system of equations which yields the required result.
		For the sake of completeness we include the calculations here.
		Write $b$ as follows (it is always possible by \red{ref}):
		\begin{equation*}
			\begin{aligned}
				b &= x_{\alpha_2}(b_2) x_{\alpha_1}(b_1) x_{\alpha_1+\alpha_2}(b_{12}) 
				x_{\alpha_{\max}-\alpha_1-\alpha_2}(b_{m12}) 
					x_{\alpha_{\max}-\alpha_1}(b_{m1}) x_{\alpha_{\max}}(b_{m}) 
					x_{\alpha_{\max}-2\alpha_1-2\alpha_2}(b_{m1122})\\
				&\quad\cdot
					x_{\alpha_{\max}-2\alpha_1-\alpha_2}(b_{m112}) x_{\alpha_{\max}-2\alpha_1}(b_{m11}) 
					\cdot \prod x_\alpha (b_\alpha),
			\end{aligned}
		\end{equation*}
		where the product is taken over all the roots in the set
		\begin{align*}	
			\Phi \setminus \{ &\alpha_1, \alpha_2, \alpha_1+\alpha_2, 
			\alpha_{\max},\alpha_{\max}-\alpha_1,\alpha_{\max}-\alpha_1-\alpha_2,\\
			&
			\alpha_{\max}-2\alpha_1, \alpha_{\max}-2\alpha_1-\alpha_2, \alpha_{\max}-2\alpha_1-2\alpha_2 \}.
		\end{align*}
		Recall that by \red{the corollary of} Lemma \ref{lemma:PC:preserves:P:i:k} we have
		\begin{equation}
			\label{eq:lemma:laundry:nonzero}
			a_1, b_2, c_{12} \in F^* \text{ and } b_1 = 0
		\end{equation}
		Let 
		\begin{equation*}
			\begin{aligned}
				\bar{b} &= x_{\alpha_2}(b_2) x_{\alpha_1+\alpha_2}(b_{12}) 
				x_{\alpha_{\max}-\alpha_1-\alpha_2}(b_{m12}) 
					x_{\alpha_{\max}-\alpha_1}(b_{m1}) x_{\alpha_{\max}}(b_{m}) \\ 
				&\quad\cdot
					x_{\alpha_{\max}-2\alpha_1-2\alpha_2}(b_{m1122})
					x_{\alpha_{\max}-2\alpha_1-\alpha_2}(b_{m112}) x_{\alpha_{\max}-2\alpha_1}(b_{m11}).
			\end{aligned}
		\end{equation*}
		It is clear that $[a,b]=[a,\bar{b}]$ and $[c,b]=[c,\bar{b}]$. 
		Now, an explicit computation shows that
		\begin{equation}
			\label{eq:lemma:laundry:1:12}
			\begin{aligned}
				e = [a,c] &= x_{\alpha_{\max}}(-2 a_1 a_{m112} c_{12} - a_{m1122} c_{12}^2 + 
					2 a_1 c_{m1} + a_1^2 c_{m11} + 2 a_1 c_{12} c_{m112})\\
				&\quad\cdot
				x_{\alpha_{\max}-\alpha_1}(-a_{m112} c_{12} + a_1 c_{m11}) \\
				&\quad\cdot
				x_{\alpha_{\max}-\alpha_1-\alpha_2}(-a_{m1122} c_{12} + a_1 c_{m112})
			\end{aligned}.
		\end{equation}
		The parameters of each transvections on the right-hand side of \eqref{eq:lemma:laundry:1:12}
		must equal zero, which gives us three equations. Our goal is to deduce an equation
		on variables $a_1,b_2,c_{12}$ and $c_{m11}$.
		Rewrite the equations obtained form \eqref{eq:lemma:laundry:1:12} as follows (form the last to the first one):
		\begin{equation}
			\label{eq:lemma:laundry:1:12:eqs}
			\begin{aligned}
				a_{m1122} c_{12} &= a_1 c_{m112} \\
				a_{m112} c_{12} &= a_1 c_{m11} \\
				0 &= - a_1^2 c_{m11} + 2 a_1 c_{m1} + a_1 c_{12} c_{m112}
				\end{aligned}.
		\end{equation}
		Next we equality
		\begin{equation}
			\label{eq:lemma:laundry:12:2}
			\begin{aligned}
				e = [c, \bar{b}] &= 
				x_{\alpha_{\max}}(2 b_{12} b_{m1122} c_{12} + 2 b_{m12} c_{12} + b_{m1122} c_{12}^2 \\
				&\qquad\qquad
				- b_{12}^2 c_{m1122} - 2 b_{12} c_{12} c_{m1122} - 2 b_{12} c_{m12})\\
				&\quad\cdot
				x_{\alpha_{\max}-\alpha_1}(b_{m112} c_{12} + b_2 b_{m1122} c_{12} - b_{12} c_{m112} \\
				&\qquad\qquad
				- b_{12} b_2 c_{m1122} - b_2 c_{12} c_{m1122} - b_2 c_{m12})\\
				&\quad\cdot
				x_{\alpha_{\max}-\alpha_1-\alpha_2}(b_{m1122} c_{12} - b_{12} c_{m1122})\\
				&\quad\cdot
				x_{\alpha_{\max}-2\alpha_1}(-2 b_2 c_{m112} - b_2^2 c_{m1122})\\
				&\quad\cdot
				x_{\alpha_{\max}-2\alpha_1-\alpha_2}(-b_2 c_{m1122})
			\end{aligned}.
		\end{equation}
		gives us five more equations. The last three equations obtained from \eqref{eq:lemma:laundry:12:2},
		as $2 \in F^*$ and $c_{12}, b_2 \in F^*$ yield
		\begin{align}
			\label{eq:lemma:laundry:cm112+bm1122+cm1122}
			c_{m1122} &= 0, &
			b_{m1122} &= 0 &&\text{ and } & c_{m112} &= 0.
		\end{align}
		This allows rewriting \eqref{eq:lemma:laundry:1:12:eqs} (keeping in mind \eqref{eq:lemma:laundry:nonzero}) as
		\begin{equation}
			\label{eq:lemma:laundry:1:12:eqs'}
			\begin{aligned}
				a_{m1122} &= 0 \\
				a_{m112} c_{12} &= a_1 c_{m11} \\
				a_1 c_{m11} &=  2 c_{m1}.
				\end{aligned}
		\end{equation}
		The second line of \eqref{eq:lemma:laundry:12:2} together with \eqref{eq:lemma:laundry:cm112+bm1122+cm1122} yields
		\begin{equation}
			\label{eq:lemma:laundry:12:2:eq2}
			b_{m112} c_{12} - b_2 c_{m12} = 0.
		\end{equation}
		Finally, the first line of \eqref{eq:lemma:laundry:12:2} together with \eqref{eq:lemma:laundry:cm112+bm1122+cm1122}
		and assumption $2 \in F^*$ yields
		\begin{equation}
			\label{eq:lemma:laundry:12:2:eq1}
			b_{m12} c_{12} - b_{12} c_{m12} = 0.
		\end{equation}
		Finally, 
		\begin{equation}
			\label{eq:lemma:laundry:1:2}
			\begin{aligned}
				e = [a, \bar{b}].c^{-1} &= 
				x_{\alpha_1+\alpha_2}(a_1 b_2 - c_{12}) \\
				&\quad\cdot
				x_{\alpha_{\max}}(-2 a_1 a_{m112} b_12 - a_{m1122} b_{12}^2 - 2 a_1^2 a_{m112} b_2 + a_1^2 a_{m1122} b_2^2\\
				&\qquad\qquad
				+ 2 a_1 b_{m1} + a_1^2 b_{m11} + 2 a_1 b_{12} b_{m112} - 2 a_{m1122} b_{12} c_{12} \\
				&\qquad\qquad 
				- 2 a_1 a_{m1122} b_2 c_{12} + 2 a_1 b_{m112} c_{12} - c_m - 2 c_{12} c_{m12})\\
				&\quad\cdot
				x_{\alpha_{\max}-\alpha_1}(-a_{m112} b_{12} - 2 a_1 a_{m112} b_2 - a_{m1122} b_{12} b_2 + a_1 b_{m11}\\
				&\qquad\qquad + a_1 b_2 b_{m112} -  a_{m1122} b_2 c_{12} - c_{m1} - c_{12} c_{m112})\\
				&\quad\cdot
				x_{\alpha_{\max}-\alpha_1-\alpha_2}(-a_{m1122} b_{12} - a_1 a_{m1122} b_2 + a_1 b_{m112} - c_{m12})\\
				&\quad\cdot
				x_{\alpha_{\max}-2\alpha_1}(-2 a_{m112} b_2 - a_{m1122} b_2^2 - c_{m11})\\
				&\quad\cdot
				x_{\alpha_{\max}-2\alpha_1-\alpha_2}(-a_{m1122} b_2 - c_{m112}).
			\end{aligned}
		\end{equation}
		The equation obtained from the last elementary root unipotent in \eqref{eq:lemma:laundry:1:2}
		doesn't give us any new information. The first one gives
		\begin{equation}
			\label{eq:lemma:laundry:1:2:c12-a1b2}
			c_{12} = a_1 b_2;
		\end{equation}
		The penultimate one together with
		\eqref{eq:lemma:laundry:1:12:eqs'} amounts to
		\begin{equation*}
			\label{eq:lemma:laundry:1:2:eq5}
			-2 a_{m112} b_2 = c_{m11}.
		\end{equation*}
		Together with \eqref{eq:lemma:laundry:nonzero} and \eqref{eq:lemma:laundry:1:2:c12-a1b2} and the assumption that $3 \in F^*$ this infers
		\begin{equation}
			\label{eq:lemma:laundry:am112}
			a_{m112} = 0.
		\end{equation}
		Finally, substituting \eqref{eq:lemma:laundry:am112} into the second line of
		\eqref{eq:lemma:laundry:1:12:eqs'} we get
		\begin{equation}
			\label{eq:lemma:laundry:cm11}
			c_{m11} = 0.
		\end{equation}
		Collecting \eqref{eq:lemma:laundry:am112}, \eqref{eq:lemma:laundry:cm112+bm1122+cm1122},
		\eqref{eq:lemma:laundry:cm11} and the first equality in \eqref{eq:lemma:laundry:1:12:eqs'}
		we get the conclusion of the lemma.

	\end{proof}
\end{lemma}

\begin{lemma}
	\label{lemma:U1}
	Assume $\rk(\Phi) \geq 3$. Let $\phi \in \PC(\Up)$. Then there exists an automorphism $\psi \in \Aut(\Up)$ 
	such that
	\begin{equation}
		\label{eq:lemma:U1:statement}
		\begin{aligned}
			\psi \circ  \phi(x_{\alpha}(1)) &\in X_{\alpha} \cdot X_{\alpha_{\max}}
		\end{aligned}
	\end{equation}
	for all $\alpha \in \Phi$ such that $m_1(\alpha) \geq 1$.
	Further, $\psi$ is the composition of an inner and an extremal automorphism.
	\begin{proof}
		First we will concentrate on the images of $x_{\alpha_1}(1)$ and $x_{\alpha_1+\alpha_2}(1)$. By
		Lemma \ref{lemma:T12} we have
		\begin{equation}
			\label{eq:lemma:U1:a1}
			a = \phi(x_{\alpha_1}(1)) \in U_1 \cdot X_{\alpha_{\max}-2\alpha_1} 
			\cdot X_{\alpha_{\max}-2\alpha_1-\alpha_2} \cdot X_{\alpha_{\max}-2\alpha_1-2\alpha_2}.
		\end{equation}
		It follows immediately that
		\begin{equation}
			\label{eq:lemma:U1:a12}
			\begin{aligned}
				b = \phi(x_{\alpha_1+\alpha_2}(1)) &= [\phi(x_{\alpha_1}(\xi)), \phi(x_{\alpha_2}(\xi))]\\
				&\le 
				[\U_1 \cdot X_{\alpha_{\max} - 2\alpha_1} 
					\cdot X_{\alpha_{\max}-2\alpha_1-\alpha_2}
					\cdot X_{\alpha_{\max}-2\alpha_1-2\alpha_2}, P^1_2] \\\
				&\le
				\U_1^{(2)} \cdot X_{\alpha_{\max}-2\alpha_1} \cdot X_{\alpha_{\max}-2\alpha_1-\alpha_2}.
			\end{aligned}
		\end{equation}
		Consider the inner automorphism $\Int_{C_1}$ defined by the matrix
		\[
			C_1 = \prod_{|j|>2} x_{\alpha_{2j}}(-(N_{\alpha_{2j},\alpha_1} a_{12})^{-1} a_{1,j'}) \cdot
			\prod_{|j|>3} x_{\alpha_{3j}}(-(N_{\alpha_{3j},\alpha_1+\alpha_2} b_{13})^{-1} b_{1,j'}),
		\]
		where $j' = \begin{cases} j & j>0 \\ 2n+1-j & j<0 \end{cases}$.
		It can be verified by a direct calculation that 
		\begin{align*}
			\Int_{C_1} \circ \phi(x_{\alpha_1}(1)) &\in X_{\alpha_1} \cdot X_{\alpha_{\max}} \cdot X_{\alpha_{\max}-\alpha_1} 
			\cdot X_{\alpha_{\max}-2\alpha_1} \cdot X_{\alpha_{\max}-2\alpha_1-\alpha_2} 
			\cdot X_{\alpha_{\max}-2\alpha_1-2\alpha_2}\\
			\Int_{C_1} \circ  \phi(x_{\alpha_1+\alpha_2}(1)) &\in X_{\alpha_1+\alpha_2} 
			\cdot X_{\alpha_{\max}} \cdot X_{\alpha_{\max}-\alpha_1} 
			\cdot X_{\alpha_{\max}-\alpha_1-\alpha_2}
			\cdot X_{\alpha_{\max}-2\alpha_1} \cdot X_{\alpha_{\max}-2\alpha_1-\alpha_2}.
		\end{align*}
		and thus by Lemma \ref{lemma:laundry} also
		\begin{align*}
			\Int_{C_1} \circ \phi(x_{\alpha_1}(1)) &\in X_{\alpha_1} \cdot X_{\alpha_{\max}} \cdot X_{\alpha_{\max}-\alpha_1} 
			\cdot X_{\alpha_{\max}-2\alpha_1} \\
			\Int_{C_1} \circ  \phi(x_{\alpha_1+\alpha_2}(1)) &\in X_{\alpha_1+\alpha_2} 
			\cdot X_{\alpha_{\max}} \cdot X_{\alpha_{\max}-\alpha_1} 
			\cdot X_{\alpha_{\max}-\alpha_1-\alpha_2}.
		\end{align*}
		It is clear that there exists an extremal automorphism $\varphi$ of the second kind such that 
		\[
			\varphi \circ \Int_{C_1} (x_{\alpha_1}(1)) \in X_{\alpha_1} \cdot X_{\alpha_{\max}} 
			\cdot X_{\alpha_{\max}-\alpha_1}.
		\]
		Let $b = \varphi \circ \Int_{C_1} (x_{\alpha_1}(1))$. Set
		\[
			C_2 = x_{\alpha_{3,-3}}(-(N_{\alpha_1+\alpha_2,\alpha_{3,-3},1} b_{13})^{-1} b_{1,-3})
		\] and $\phi_2 = \Int_{C_2} \circ \phi \circ \Int_{C_1}$. It is easily seen that
		\begin{align*}
			\phi_2(x_{\alpha_1}(1)) &\in X_{\alpha_1} \cdot X_{\alpha_{\max}} \cdot X_{\alpha_{\max}-\alpha_1} \\
			\phi_2(x_{\alpha_1+\alpha_2}(1)) &\in X_{\alpha_1+\alpha_2} 
			\cdot X_{\alpha_{\max}} \cdot X_{\alpha_{\max}-\alpha_1}.
		\end{align*}
		Finally, by computing the commutator $[\phi_2(x_{\alpha_1}(1)),\phi_2(x_{\alpha_1+\alpha_2}(1))]$, which
		must equal $e$ we deduce that $\phi_2(x_{\alpha_1+\alpha_2}(1)) \in X_{\alpha_1+\alpha_2} \cdot X_{\alpha_{\max}}$.
		
		We continue constructing the required automorphism $\psi_k$ by induction. Let $k < n$ and assume that
		\begin{align*}
			\phi_{k-1}(x_{\alpha_1}(1)) &\in X_{\alpha_1} \cdot X_{\alpha_{\max}} \cdot X_{\alpha_{\max}-\alpha_1} \\
			\phi_{k-1}(x_{\alpha_1+\dots+\alpha_j}(1)) &\in X_{\alpha_1+\dots+\alpha_j} 
			\cdot X_{\alpha_{\max}}.
		\end{align*}
		for all $j<k$. We will show that the same holds for $k=k+1$. Let $b = \psi_{k-1}(x_{\alpha_1+\dots+\alpha_k}(1))$.
		Set
		\[
			\bar{C}_k = \prod_{|j|>k+1} x_{\alpha_{k+1,j}}(-(N_{\alpha_1+\dots+\alpha_k, \alpha_{k+1,j}} 
			b_{1,k+1})^{-1}) b_{1,j'}),
		\]
		where, as usual, $j' = \begin{cases} j & j>0 \\ 2n+1-j & j<0 \end{cases}$. Now let
		$c = \Int_{\bar{C}_k} \circ \psi_{k-1} (x_{\alpha_1+\dots+\alpha_k}(1)).$ Set
		\[
			C_k = x_{k+1,-(k+1)}(-(N_{\alpha_1+\dots+\alpha_k,\alpha_{k+1,-(k+1)},1} c_{1,k+1})^{-1}) c_{1,2n+1-(k+1)}) 
			\cdot \bar{C}_k.
		\]
		Set $\psi_k = \Int_{C_k} \circ \psi_{k-1}$ and $d = \psi_k(x_{\alpha_1+\dots+\alpha_k}(1))$. It is easily
		checked that $c_{1,k+2}=\dots=c_{1,2n+1-(k+1)}=0$. Further, as $C_k$ centralizes the subgroups
		\[
			X_{\alpha_1} \cdot X_{\alpha_{\max}} \cdot X_{\alpha_{\max}-\alpha_1} \text{ and }
			X_{\alpha_1+\dots+\alpha_j} \cdot X_{\alpha_{\max}}
		\]
		for all $j<k$ it follows that $\psi_k(x_{\alpha_1+\dots+\alpha_j}(1))=\psi_{k-1}(x_{\alpha_1+\dots+\alpha_j}(1))$
		for $j<k$. Finally, by considering the commutators (which much be zero) of $\psi_k(x_{\alpha_1+\dots+\alpha_k})$
		with $\psi_k(x_{\alpha_1+\dots+\alpha_j}), j<k$ we conclude that $c_{1,2n+1-(k+1)+1}=\dots=c_{1,2n-1}=0$,
		thus 		
		\begin{align*}
			\phi_{k}(x_{\alpha_1}(1)) &\in X_{\alpha_1} \cdot X_{\alpha_{\max}} \cdot X_{\alpha_{\max}-\alpha_1} \\
			\phi_{k}(x_{\alpha_1+\dots+\alpha_j}(1)) &\in X_{\alpha_1+\dots+\alpha_j} 
			\cdot X_{\alpha_{\max}}.
		\end{align*}
		for all $j \leq k$. By induction this holds for all $k$ up to $k = n-1$.
		
		Finally, by computing the (necessarily trivial) commutator of
		\[
			\phi_{n-1}(x_{\alpha_{\max}-\alpha_1-\dots-\alpha_j}(1))
		\]
		with $\phi_{n-1}(x_{\alpha_1-\dots-\alpha_{j-l}}(1))$, $1 \leq l \leq j-1$ and $1 \leq j \leq n-1$
		we deduce that
		\begin{align*}
			\phi_{n-1}(x_{\alpha_1}(1)) &\in X_{\alpha_1} \cdot X_{\alpha_{\max}} \cdot X_{\alpha_{\max}-\alpha_1} \\
			\phi_{n-1}(x_{\alpha}(1)) &\in X_{\alpha} \cdot X_{\alpha_{\max}}.
		\end{align*}
		for all $\alpha \in \Phi$ such that $m_1(\alpha) = 1$. 
		Finally, an appropriately chosen extremal automorphism $\eta$ sends $\phi_{n-1} x_{\alpha_1}(1))$ to
		$X_{\alpha_1} \cdot X_{\alpha_{\max}}$ and leaves $\phi_{n-1} x_{\alpha}(1))$, 
		$m_1(\alpha)\geq 1$, $\alpha \neq \alpha_1$ invariant. Set $\psi_n = \eta \circ \psi_{n-1}$.
		We have shown that
		\begin{align*}
			\phi_{n}(x_{\alpha}(1)) &\in X_{\alpha} \cdot X_{\alpha_{\max}}
		\end{align*}
		for all $\alpha \in \Phi$ such that $m_1(\alpha)\geq 1$.
		Notice that $\psi_{n} = \eta \circ \Int_{C_{n-1} \cdots C_{2}} \circ \varphi \circ \Int_{C_1} \circ \phi$.
		As the subgroup of internal automorphisms normalized by the subgroup of extremal automorphisms
		it follows that $\psi = \psi_{n} \circ \phi^{-1}$ is, as required, a composition of an internal and 
		an extremal automorphism. This completes the proof.
	\end{proof}			
\end{lemma}

\begin{lemma}
	\label{lemma:simple:roots}
	Assume $\rk(\Phi) \geq 3$ and $2F=F$. Let $\phi \in \PC(\Up)$. Suppose that
	\begin{equation}
		\label{eq:lemma:simple:roots:assumption}
		\begin{aligned}
			\phi(x_{\alpha}(1)) &\in X_{\alpha} \cdot X_{\alpha_{\max}}
		\end{aligned}
	\end{equation}
	for all $\alpha\in \Phi$ such that $m_1(\alpha) \geq 1$.
	Then there exist such an inner automorphism $\Int_C$ such that
	\begin{equation}
		\label{eq:lemma:simple:roots:conclusion1}
		\begin{aligned}
			\Int_C \circ \phi(X_{\alpha}) \subseteq X_{\alpha} \cdot X_{\alpha_{\max}}
		\end{aligned}
	\end{equation}	
	for all short simple roots $\alpha \in \{\alpha_1,\dots,\alpha_{n-1}\}$ as well as for $\alpha = \alpha_{n-1}+\alpha_n$.
	Further,
	\begin{equation}
		\label{eq:lemma:simple:roots:conclusion2}
		\begin{aligned}
			\Int_C \circ \phi(X_{\alpha_{1i}}) \subseteq X_{\alpha_{1i}}
		\end{aligned}
	\end{equation}	
	for all $2 \leq i \leq 2n-2$. Finally, 
	\[
		\Int_C \circ \phi(X_{\alpha_n}) \subseteq X_{\alpha_n} \cdot X_{\alpha_1+\dots+\alpha_n}
		\cdot X_{\alpha_{\max}}.
	\]
	\begin{proof}
		Let $1 \leq j \leq n$. By assumption, $\phi$ maps $x_{\alpha_1+\dots+\alpha_i}(1)$ to 
		$X_{\alpha_1+\dots+\alpha_i} \cdot X_{\alpha_{\max}}$. Combining Lemma \ref{lemma:simple:roots:as:C}
		and Corollary \ref{cor:centralizers} we conclude that the subgroups
		\[
			X_{\alpha_i} \cdot X_{\alpha_1+\dots+\alpha_i} \cdot 
				X_{\alpha_{\max}-\alpha_1-\dots-\alpha_{i-1}} \cdot X_{\alpha_{\max}}
		\]
		are preserved by $\phi$ for all $i \leq n$. In particular,
		\begin{align*}
			\phi(x_{\alpha_1}(1)) &= x_{\alpha_1}(a_1)
				x_{\alpha_{\max}}(*) \\
			\phi(x_{\alpha_i}(1)) &= x_{\alpha_i}(a_i) x_{\alpha_1+\dots+\alpha_i}(b_i) 
				x_{\alpha_{\max}-\alpha_1-\dots-\alpha_{i-1}}(c_i) x_{\alpha_{\max}}(*), \qquad 2 \leq i \leq n-1 \\
		\end{align*}
		Let 
		\begin{align*}
			C_1 &= \prod_{i=2}^{n-1} x_{\alpha_1+\dots+\alpha_{i-1}}	
				(-(N_{\alpha_1+\dots+\alpha_{i-1}, \alpha_i} a_i)^{-1} b_i)\\
			&\quad\cdot 
			\prod_{i=2}^{n-1} x_{\alpha_{\max}-\alpha_1-\dots-\alpha_i}
				(-(N_{\alpha_{\max}-\alpha_1-\dots-\alpha_i, \alpha_i} a_i)^{-1} c_i).
		\end{align*}
		It is easily seen that
		\begin{align*}
			\Int_{C_1} \circ \phi(x_{\alpha_i}(1)) &= x_{\alpha_i}(a_i) x_{\alpha_{\max}}(*), \qquad 1 \leq i \leq n-1.
		\end{align*}
		Further, as $C_1$ normalizes the subgroups $X_\alpha \cdot X_{\alpha_{\max}}$ for all $\alpha \in \Phi$
		such that $m_1(\alpha) \geq 1$ it follows that the inclusions \eqref{eq:lemma:simple:roots:assumption}
		still hold for $\phi = \Int_{C_1} \circ \phi$. 
		
		Now, we can again use the combination of Lemma \ref{lemma:simple:roots:as:C} and Corollary \ref{cor:centralizers}
		in order to deduce that the subgroup
		\begin{equation*}
			X_{\alpha_n} \cdot X_{\alpha_{n-1}+\alpha_n} \cdot X_{2\alpha_{n-1}+\alpha_n} \cdot 
			X_{\alpha_1+\dots+\alpha_n} \cdot X_{\alpha_1+\dots+\alpha_n +\alpha_{n-1}} X_{\alpha_{\max}}
		\end{equation*}
		is preserved by $\Int_{C_1} \circ \phi$.
		By Lemma \ref{lemma:PC:preserves:Up:s} any $\PC$-map preserves $\Up^{(s)}$ for all $s$, thus
		the subgroup 
		\begin{equation*}
			X_{\alpha_{n-1}+\alpha_n} \cdot X_{2\alpha_{n-1}+\alpha_n} \cdot 
			X_{\alpha_1+\dots+\alpha_n} \cdot X_{\alpha_1+\dots+\alpha_n +\alpha_{n-1}} X_{\alpha_{\max}}
		\end{equation*}
		is also preserved by $\Int_{C_1} \circ \phi$. In particular,		
		\begin{align*}
			\Int_{C_1} \circ \phi(x_{\alpha_{n-1}+\alpha_n}(1)) &= x_{\alpha_{n-1}+\alpha_n}(a_n)
			x_{2\alpha_{n-1}+\alpha_n}(d_n) x_{\alpha_1+\dots+\alpha_n}(c_n)\\
			&\quad\cdot
			x_{\alpha_1+\dots+\alpha_n +\alpha_{n-1}}(b_n) x_{\alpha_{\max}}(*).
		\end{align*}
		Let
		\[
			\Int_{C_1} \circ \phi(x_{\alpha_1+\dots+\alpha_{n-2}}(1)) = 
			x_{\alpha_1+\dots+\alpha_{n-2}}(\zeta) x_{\alpha_{\max}}(*).
		\]
		By \red{ the corollary of } Lemma \ref{lemma:PC:preserves:P:i:k}, $\zeta \neq 0$.
		Now, observe that \red{ here we use an explicit computation of a structure constant! }
		\begin{equation}
			\begin{aligned}
				\Int_{C_1} \circ \phi(x_{\alpha_1+\dots+\alpha_n}(1))
				&=
				\Int_{C_1} \circ \phi([x_{\alpha_1+\dots+\alpha_{n-2}}(1), x_{\alpha_{n-1}+\alpha_n}(1)]) \\
				&= 
				[x_{\alpha_1+\dots+\alpha_{n-2}}(\zeta) x_{\alpha_{\max}}(*), 
				x_{\alpha_{n-1}+\alpha_n}(a_n)
				x_{2\alpha_{n-1}+\alpha_n}(d_n) \\
				&\qquad\cdot 
				x_{\alpha_1+\dots+\alpha_n}(c_n)
				x_{\alpha_1+\dots+\alpha_n +\alpha_{n-1}}(b_n) x_{\alpha_{\max}}(*)] \\
				&=
				x_{\alpha_1+\dots+\alpha_n}(*) x_{\alpha_1+\dots+\alpha_n+\alpha_{n-1}}(\pm \zeta d_n) x_{\alpha_{\max}}.
			\end{aligned}
		\end{equation}
		On the other hand, 
		\[
			\Int_{C_1} \circ \phi(x_{\alpha_1+\dots+\alpha_n}(1)) \in X_{\alpha_1+\dots+\alpha_n} \cdot X_{\alpha_{\max}}.
		\]
		Thus $d_n = 0$.
		Set \begin{align*}
			C_2 &= 
			x_{\alpha_1+\dots+\alpha_{n-1}}
				(-(N_{\alpha_1+\dots+\alpha_{n-1}, \alpha_{n-1}+\alpha_n} a_n)^{-1} b_n).
		\end{align*}
		It is easily checked that 
		\begin{align*}
			\Int_{C_2 C_1} \circ \phi(x_{\alpha_{n-1}+\alpha_n}(1)) &= x_{\alpha_{n-1}+\alpha_n}(a_n)
			x_{\alpha_{1}+\dots+\alpha_n}(c_n) x_{\alpha_{\max}}(*).
		\end{align*}
		Moreover, as $C_2$ normalizes $X_{\alpha} \cdot X_{\alpha_{\max}}$ whenever $m_1(\alpha) \geq 1$ or 
		$\alpha \in \Pi \setminus \{\alpha_n\}$, it follows that
		\begin{align*}
			\Int_{C_2 C_1} \circ \phi (x_{\alpha}(1)) &\in X_{\alpha} \cdot X_{\alpha_{\max}}
		\end{align*}
		for all $\alpha \in \{\beta \in \Phi \mid m_1(\beta) \geq 1 \} \cup 
		\{\alpha_1, \dots, \alpha_{n-1}\}$. 
		Replace $\phi$ with $\Int_{C_2 C_1} \circ \phi$.
		
		The next step is to show that $\phi(X_{\alpha_1}) \subseteq X_{\alpha_1} \cdot X_{\alpha_{\max}}$
		and $\phi(X_{\alpha}) \subseteq X_{\alpha}$ for all $\alpha \in \Phi$ such that $m_1(\alpha) \geq 1,
		\alpha \neq \alpha_1$. Indeed, we have already shown that the subgroups
		\[
			X_{\alpha_i} \cdot X_{\alpha_1+\dots+\alpha_i} \cdot 
				X_{\alpha_{\max}-\alpha_1-\dots-\alpha_{i-1}} \cdot X_{\alpha_{\max}}
		\]
		are preserved by $\phi$ for $1 \leq i \leq n$. In particular,
		\[
			\phi(x_{\alpha_1}(\xi)) \in \phi(X_{\alpha_1} \cdot X_{\alpha_{\max}}) \subseteq 
			X_{\alpha_1} \cdot X_{\alpha_{\max}}.
		\]
		We continue by induction on height of $\alpha$. Suppose $\phi(X_{\alpha}) \subseteq
		X_{\alpha} X_{\alpha_{\max}}$ for all $\alpha \in \Phi$ such that
		$m_1(\alpha) \geq 1$ and $\height(\alpha)<k<2n-1$. Then there exists
		a unique simple root $\alpha_j$ such that $\alpha+\alpha_j$ is a root. We will show
		that $\phi(X_{\alpha+\alpha_j}) \subseteq X_{\alpha+\alpha_j}$. Indeed, if
		$j \neq n$ then
		\[
			\phi(x_{\alpha+\alpha_j}(\xi))=[\phi(x_{\alpha}(N_{\alpha,\alpha_j}^{-1} \xi)), \phi(x_{\alpha_j}(1))]
			\in [X_{\alpha} \cdot X_{\alpha_{\max}}, X_{\alpha_j} \cdot X_{\alpha_{\max}}]
			\subseteq X_{\alpha+\alpha_j}.
		\]
		Finally, if $j = n$, then $\alpha = \alpha_1+\dots+\alpha_{n-1}$. As $n \geq 3$ we have
		$\alpha+\alpha_j = (\alpha_1+\dots+\alpha_{n-2})+(\alpha_{n-1}+\alpha_n)$. Thus
		\begin{align*}
			\phi(x_{\alpha+\alpha_j}(\xi))&
			=[\phi(x_{\alpha_1+\dots+\alpha_{n-2}}(N_{\alpha_1+\dots+\alpha_{n-2},\alpha_{n-1}+\alpha_n}^{-1} \xi)),
			\phi(x_{\alpha_{n-1}+\alpha_n}(1))] \\
			&\in 
			[X_{\alpha_1+\dots+\alpha_{n-2}} \cdot X_{\alpha_{\max}}, X_{\alpha_{n-1}+\alpha_n}
			\cdot X_{\alpha_{1}+\dots+\alpha_n} \cdot X_{\alpha_{\max}}]
			\subseteq X_{\alpha+\alpha_j}.
		\end{align*}
		Summing up, we get $\phi(X_{\alpha_1}) \subseteq X_{\alpha_1} \cdot X_{\alpha_{\max}}$
		and $\phi(X_{\alpha}) \subseteq X_{\alpha}$ for all $\alpha \in \Phi$ such that $m_1(\alpha) \geq 1,
		\alpha \neq \alpha_1$. 

		Now, observe that for some $\zeta \neq 0$:
		\begin{align*}
			\phi(x_{\alpha_1+\dots+\alpha_{n-1}+\alpha_n+\alpha_{n-1}}(1)) &= 
			[\phi(x_{\alpha_1+\dots+\alpha_{n-1}}(N^{-1}_{\alpha_1+\dots+\alpha_{n-1},\alpha_n+\alpha_{n-1}})), 
			\phi(x_{\alpha_n+\alpha_{n-1}}(1))]\\
			&=
			[x_{\alpha_1+\dots+\alpha_{n-1}}(\zeta), x_{\alpha_{n-1}+\alpha_n}(a_n)
			x_{\alpha_{1}+\dots+\alpha_n}(c_n) x_{\alpha_{\max}}(*)] \\
			&= x_{\alpha_1+\dots+\alpha_{n-1}+\alpha_n+\alpha_{n-1}}(*)
			x_{\alpha_{\max}}(\pm 2 c_n \zeta).
		\end{align*}
		On the other hand, we have just shown that 
		$\phi(x_{\alpha_1+\dots+\alpha_{n-1}+\alpha_n+\alpha_{n-1}}(1)) \in 
		X_{\alpha_1+\dots+\alpha_{n-1}+\alpha_n+\alpha_{n-1}}$. Thus $c_n = 0$ and
		\[
			\phi(x_{\alpha_{n-1}+\alpha_n}(1)) \in X_{\alpha_{n-1}+\alpha_n} \cdot X_{\alpha_{\max}}.
		\]
			
		\red{In the following step structure some of constants are computed }
		Now fix $2 \le i \leq n-1$. We will show that $\phi(X_{\alpha_i}) \subseteq X_{\alpha_i} \cdot X_{\alpha_{\max}}$.
		We already know that 
		\[
			\phi(x_{\alpha_i}(\xi)) = x_{\alpha_i}(a) x_{\alpha_1+\dots+\alpha_i}(b) x_{\alpha_{\max}-\alpha_1-\dots-\alpha_{i-1}}(c) x_{\alpha_{\max}}(*)
		\]
		for some $a,b,c \in F$. Let $\phi(x_{\alpha_1+\dots+\alpha_{i-1}}(1)) = x_{\alpha_1+\dots+\alpha_{i-1}}(d)$. Then
		\begin{align*}
			X_{\alpha_1+\dots+\alpha_i} &= \phi(X_{\alpha_1+\dots+\alpha_i}) \ni \phi(x_{\alpha_1+\dots+\alpha_i}(\xi)) = 
			\phi([x_{\alpha_1+\dots+\alpha_{i-1}}(1),x_{\alpha_i}(\xi)])\\
			& =
			[x_{\alpha_1+\dots+\alpha_{i-1}}(d) x_{\alpha_{\max}}(*), x_{\alpha_i}(a) x_{\alpha_1+\dots+\alpha_i}(b) x_{\alpha_{\max}-\alpha_1-\dots-\alpha_{i-1}}(c) x_{\alpha_{\max}}(*)]\\
			&=
			X_{\alpha_1+\dots+\alpha_i}(*)x_{\alpha_{\max}}(cd).
		\end{align*}
		Thus $cd=0$ and as $d\neq 0$ then $c=0$.
		In the same way, let $\phi(x_{\alpha_{\max}-\alpha_1-\dots-\alpha_i}(1)) = 
		x_{\alpha_{\max}-\alpha_1-\dots-\alpha_i}(f)$. Then
		\begin{align*}
			X_{\alpha_{\max}-\alpha_1-\dots-\alpha_{i-1}} &= \phi(X_{\alpha_{\max}-\alpha_1-\dots-\alpha_{i-1}}) \ni \phi(x_{\alpha_{\max}-\alpha_1-\dots-\alpha_{i-1}}(\xi)) \\
			& = 
			\phi([x_{\alpha_{\max}-\alpha_1-\dots-\alpha_{i}}(1),x_{\alpha_i}(\xi)])\\
			& =
			[x_{\alpha_{\max}-\alpha_1-\dots-\alpha_i}(f)
		x_{\alpha_{\max}}(*), x_{\alpha_i}(a) x_{\alpha_1+\dots+\alpha_i}(b)  x_{\alpha_{\max}}(*)]\\
			&=
			X_{\alpha_{\max}-\alpha_1-\dots-\alpha_{i-1}}(*)x_{\alpha_{\max}}(bf).
		\end{align*}		
		Therefore $b=0$ and $\phi(X_{\alpha_i}) \subseteq X_{\alpha_i} \cdot X_{\alpha_{\max}}$ for all
		$i \neq n$. 
		
		Finally, let
		\[
			\phi(x_{\alpha_{n-1}+\alpha_n}(\xi)) = 
				x_{\alpha_{n-1}+\alpha_n}(a) \cdot 
				x_{\alpha_n+2\alpha_{n-1}}(b) 
				x_{\alpha_1+\dots+\alpha_n}(c)
				x_{\alpha_1+\dots+\alpha_n +\alpha_{n-1}}(d)
				x_{\alpha_{\max}}(*).
		\]
		Then for some $\zeta \neq 0$:
		\begin{align*}
			X_{\alpha_{1}+\dots+\alpha_n} &\ni \phi(x_{\alpha_{1}+\dots+\alpha_n}(\xi))
			= 
			\phi([x_{\alpha_{1}+\dots+\alpha_{n-2}}(1), x_{\alpha_{n-1}+\alpha_n}(\xi))\\
			&=
			[x_{\alpha_{1}+\dots+\alpha_{n-2}}(\zeta), x_{\alpha_{n-1}+\alpha_n}(a) \cdot 
				x_{\alpha_n+2\alpha_{n-1}}(b) 
				x_{\alpha_1+\dots+\alpha_n}(c)
				x_{\alpha_1+\dots+\alpha_n +\alpha_{n-1}}(d)
				x_{\alpha_{\max}}(*)]\\
			&=
			x_{\alpha_1+\dots+\alpha_n}(*) x_{\alpha_1+\dots+\alpha_n+\alpha_{n-1}}(\pm \zeta b) x_{\alpha_{\max}}(*).
		\end{align*}
		Thus $b=0$. In the same way
		\begin{align*}
			X_{\alpha_{1}+\dots+\alpha_n+\alpha_{n-1}} &\ni \phi(x_{\alpha_{1}+\dots+\alpha_n+\alpha_{n-1}}(\xi))
			= 
			\phi([x_{\alpha_{1}+\dots+\alpha_{n-1}}(1), x_{\alpha_{n-1}+\alpha_n}(\xi))\\
			&=
			[{\alpha_{1}+\dots+\alpha_{n-1}}(\zeta), x_{\alpha_{n-1}+\alpha_n}(a) \cdot 
				x_{\alpha_1+\dots+\alpha_n}(c)
				x_{\alpha_1+\dots+\alpha_n +\alpha_{n-1}}(d)
				x_{\alpha_{\max}}(*)]\\
			&=
			x_{\alpha_1+\dots+\alpha_n+\alpha_{n-1}}(*) x_{\alpha_{\max}}(2 \zeta c).
		\end{align*}
		Thus $c=0$. 
		Finally, again consider the commutator
		\begin{align*}
			X_{\alpha_{1}+\dots+\alpha_n} &\ni \phi(x_{\alpha_{1}+\dots+\alpha_n}(\xi))
			= 
			\phi([x_{\alpha_{1}+\dots+\alpha_{n-2}}(1), x_{\alpha_{n-1}+\alpha_n}(\xi))\\
			&=
			[x_{\alpha_{1}+\dots+\alpha_{n-2}}(\zeta), x_{\alpha_{n-1}+\alpha_n}(a) \cdot 
				x_{\alpha_1+\dots+\alpha_n +\alpha_{n-1}}(d)
				x_{\alpha_{\max}}(*)]\\
			&=
			x_{\alpha_1+\dots+\alpha_n}(*) x_{\alpha_{\max}}(\pm 2 \zeta d).
		\end{align*}
		Thus, $d = 0$ and $\phi(X_{\alpha_{n-1}+\alpha_n}) \subseteq X_{\alpha_{n-1}+\alpha_n} \cdot X_{\alpha_{\max}}$.
		This completes the proof.
	\end{proof}
\end{lemma}

\begin{lemma}
	\label{lemma:X:dirty}
	Assume $\rk(\Phi) \geq 3$. Let $\phi \in \PC(\Up)$. Suppose that
	\begin{equation}
		\label{eq:lemma:X:dirty:assumption}
		\begin{aligned}
			\phi(X_{\alpha}) \in X_{\alpha} \cdot X_{\alpha_{\max}}
		\end{aligned}
	\end{equation}
	for all $\alpha \in \Pi' = \{\alpha_1,\dots,\alpha_{n-1},\alpha_{n-1}+\alpha_n\}$.
	Then there exists a diagonal automorphism $D_1$ and a quasi-diagonal automorphism $D_2$ such that
	\begin{equation}
		\label{eq:lemma:X:dirty:conclusion:simple}
		\begin{aligned}
			D_2 \circ D_1 \circ \phi(x_{\alpha}(1)) = x_\alpha(1) \cdot x_{\alpha_{\max}}(*),
		\end{aligned}
	\end{equation}
	for all $\alpha \in \Pi'$;
	\begin{equation}
		\label{eq:lemma:X:dirty:conclusion:short}
		\begin{aligned}
			D_2 \circ D_1 \circ \phi(x_{\alpha}(1)) = x_\alpha(1)
		\end{aligned}
	\end{equation}
	for all $\alpha \in \Phi^S \setminus \Pi'$ and
	\begin{equation}
		\label{eq:lemma:X:dirty:conclusion:long}
		\begin{aligned}
			D_2 \circ D_1 \circ \phi(x_{\alpha}(2)) = x_\alpha(2)
		\end{aligned}
	\end{equation}
	for all $\alpha \in \Phi^L \setminus \{\alpha_n\}$.	
	\begin{proof}
		By assumption,
		\[
			\phi(x_{\alpha_i}(1)) = x_{\alpha_i}(d_i) \cdot x_{\alpha_{\max}}(*)
		\]
		for all $1 \leq i \leq n-1$.
		Clearly, the diagonal automorphism $D_1 = h_1(1) h_2(d_1) h_3(d_1 d_2) \cdot \dots \cdot h_n(\prod_{i=1}^{n-1} d_i)$
		satisfies the following condition:
		\begin{equation}
			\label{eq:lemma:X:dirty:1}
			D_1 \circ \phi(x_{\alpha_i}(1)) = x_{\alpha_i}(1) \cdot x_{\alpha_{\max}}(*)
		\end{equation}
		for all $1 \leq i \leq n-1$. Now,
		\[
			\phi(x_{\alpha_{n-1}+\alpha_n}(1)) = x_{\alpha_{n-1}+\alpha_n}(d_n) \cdot x_{\alpha_{\max}}(*).
		\]
		Set $D_2$ be the quasi-diagonal automorphism induced by $d_n^{-1}$.
		Clearly, $D_1$ and $D_2$ satisfy \eqref{eq:lemma:X:dirty:conclusion:simple} for $\alpha$ being a simple short root
		and the root $\alpha_{n-1}+\alpha_n$. 

		The inclusions \eqref{eq:lemma:X:dirty:conclusion:short} and
		\eqref{eq:lemma:X:dirty:conclusion:long} are obtained by induction on height.
		Indeed, every root in $\Phi \setminus (\Pi \cup \{\alpha_{n-1}+\alpha_n\})$ is a sum of two short roots.
		Indeed, let $\alpha_{ij} \in \Phi \setminus (\Pi \cup \{\alpha_{n-1}+\alpha_n\})$. 
		If $j>0$, then, as $\height(\alpha_{ij}) \geq 2$, $\alpha_{ij} = \alpha_i+\alpha_{i+1,j}$.
		Assume $j<0$. Then $\alpha_{ij}$ can be of height 2 only if $\alpha_{ij} = \alpha_{n-1}+\alpha_n$.
		Thus $\height(\alpha_{ij}) \geq 3$. It follows that $i \leq n-1$.
		If $j>i+1$, then $\alpha_{ij} = \alpha_i + \alpha_{i+1,j}$ and $\alpha_{i+1,j}$ is short.
		Finally, if $j = i+1$ then $\alpha_{ij} = \alpha_{i,i+2}+\alpha_{i+2,j}$ and
		both $\alpha_{i,i+2}$ and $\alpha_{i+2,j}$ are short.
		\red{ uses explicit structure constants:}
		Therefore,
		\[
			\phi(x_{\alpha_{ij}}(C)) = \phi([x_{\alpha_{i,i'}}(1),x_{\alpha_{i',j}}(1)]) =
			[x_{\alpha_{i,i'}}(1),x_{\alpha_{i',j}}(1)] = x_{\alpha_{ij}}(C),
		\]
		where $C = 1$, if $\alpha_{ij}$ is short and $C=2$, otherwise, and $i' \in \{i+1,i+2\}$
		is chosen as above.
	\end{proof}
\end{lemma}

\red{Proper introduction} The next lemma is actually one of the most difficult parts of the proof.
We suggest the reader to go through Lemma \ref{lemma:field:autom} first and come back to Lemma
\ref{lemma:n:1:to:n} afterwards...

\begin{lemma}
	\label{lemma:n:1:to:n}
	Assume $n \geq 4$ \red{4 is important here, see comment below} and $2F=F$. Let $\phi \in \PC(\Up)$. Suppose that
	\begin{align*}
		\phi(x_{\alpha_n}(\xi)) &= x_{\alpha_n}(\xi) \cdot x_{\alpha_1+\dots+\alpha_n}(*) \\
		\phi(x_{\alpha_{n-1}+\alpha_n}(\xi)) &= x_{\alpha_{n-1}+\alpha_n}(\xi) \cdot x_{\alpha_{\max}}(*) \\			
		\phi(x_{\alpha}(\xi)) &= x_{\alpha}(\xi)
	\end{align*}
	for all $\alpha \neq \alpha_n, \alpha_{n-1}+\alpha_n$ and all $\xi \in F$.
	Then
	\begin{align*}
		\phi(x_{\alpha_n}(\xi) x_{\alpha_1+\dots+\alpha_n}(\zeta)) &= x_{\alpha_n}(\xi) x_{\alpha_1+\dots+\alpha_n}(\zeta)
		x_{\alpha_{\max}}(*), \\
		\phi(x_{\alpha_1+\dots+\alpha_n}(\xi)x_{\alpha_{\max}}(\zeta))
		&=
		x_{\alpha_1+\dots+\alpha_n}(\xi)x_{\alpha_{\max}}(\zeta)			
	\end{align*}
	for all $\xi,\zeta \in F$.
	\begin{proof}
		The combination of Lemma \ref{lemma:simple:roots:as:C} and Corollary \ref{cor:centralizers}
		yields that
		\[
				\phi(X_{\alpha_n} \cdot X_{\alpha_1+\dots+\alpha_n} \cdot X_{\alpha_{\max}})
				\subseteq X_{\alpha_n} \cdot X_{\alpha_1+\dots+\alpha_n} \cdot X_{\alpha_{\max}}.
		\] 
		In particular,
		\begin{align*}
			\phi(x_{\alpha_n}(\xi) x_{\alpha_1+\dots+\alpha_n}(\zeta)) =
			x_{\alpha_n}(a) x_{\alpha_1+\dots+\alpha_n}(b) x_{\alpha_{\max}}(*)
		\end{align*}
		for some $a,b \in F$. By Lemma \ref{lemma:PC:preserves:P:i:k}, $a \neq 0$.
		Next, observe that
		\begin{align*}
			x_{\alpha_{\max}}(&N_{\alpha_1+\dots+\alpha_{n-1}, \alpha_1+\dots+\alpha_n} 
				N_{\alpha_1+\dots+\alpha_{n-1}, \alpha_n,1} \xi)\\
			&=
			\phi(x_{\alpha_{\max}}(N_{\alpha_1+\dots+\alpha_{n-1}, \alpha_1+\dots+\alpha_n} 
				N_{\alpha_1+\dots+\alpha_{n-1}, \alpha_n,1} \xi))\\
			&= 
			\phi([x_{\alpha_1+\dots+\alpha_{n-1}}(1), [x_{\alpha_1+\dots+\alpha_{n-1}}(1), 
				x_{\alpha_n}(\xi)x_{\alpha_1+\dots+\alpha_n}(\zeta))]])\\
			&=
			[\phi(x_{\alpha_1+\dots+\alpha_{n-1}}(1)), 
				[\phi(x_{\alpha_1+\dots+\alpha_{n-1}}(1)), \phi(x_{\alpha_n}(\xi)x_{\alpha_1+\dots+\alpha_n}(\zeta))]] \\
			&= [x_{\alpha_1+\dots+\alpha_{n-1}}(1), 
				[x_{\alpha_1+\dots+\alpha_{n-1}}(1), x_{\alpha_n}(a) x_{\alpha_1+\dots+\alpha_n}(b) x_{\alpha_{\max}}(*)]]\\
			&= [(x_{\alpha_1+\dots+\alpha_{n-1}}(1), 
				x_{\alpha_1+\dots+\alpha_n}(N_{\alpha_1+\dots+\alpha_{n-1}, \alpha_n,1} a) x_{\alpha_{\max}}(*)]\\
			&=  x_{\alpha_{\max}}(N_{\alpha_1+\dots+\alpha_{n-1}, \alpha_1+\dots+\alpha_n} 
				N_{\alpha_1+\dots+\alpha_{n-1}, \alpha_n,1} a).
		\end{align*}
		Therefore, $a = \xi$. 
		
		The proof that $b = \zeta$ is quite some more complicated. The problem is, we can't \textit{extract}
		an elementary root unipotent with parameter $\zeta$ by taking commutators of
		$\phi(x_{\alpha_n}(\xi)x_{\alpha_1+\dots+\alpha_n}(\zeta))$ with elementary root unipotents.
		To aid this problem, we use another method: rather than extracting elements, we will
		construct a certain auxiliary matrix, as a commutator of matrices which are known to be preserved by $\phi$.
		In turn, this matrix contains enough info to infer the equality $b=\zeta$. \red{Rewrite!!!!!!}
		
		First, we will show that for all $\mu,\nu \in F$
		\begin{equation}
			\label{eq:lemma:n:1:to:n:a1+a12}
			\phi(x_{\alpha_1}(\mu) x_{\alpha_1+\alpha_2}(\nu)) = 
				x_{\alpha_1}(\mu) x_{\alpha_1+\alpha_2}(\nu) x_{\alpha_{\max}}(*).
		\end{equation}
		Indeed, by Lemma \ref{lemma:simple:roots:as:C}, the subgroup
		\[
			X_{\alpha_1} \cdot 	X_{\alpha_1+\alpha_2} \cdot X_{\alpha_{\max}}		
		\]
		is preserved by $\phi$. In particular,
		\[
			\phi(x_{\alpha_1}(\mu) x_{\alpha_1+\alpha_2}(\nu)) = 
				x_{\alpha_1}(u) x_{\alpha_1+\alpha_2}(v) x_{\alpha_{\max}}(*).
		\]
		Then
		\begin{align*}
			x_{\alpha_1+\alpha_2}(N_{\alpha_1,\alpha_2}\mu)		
			&=
			\phi(x_{\alpha_1+\alpha_2}(N_{\alpha_1,\alpha_2}\mu))
			=
			[\phi(x_{\alpha_1}(\mu) x_{\alpha_1+\alpha_2}(\nu)), \phi(x_{\alpha_2}(1))]\\
			&=
			[x_{\alpha_1}(u) x_{\alpha_1+\alpha_2}(v) x_{\alpha_{\max}}(*), x_{\alpha_2}(1)]
			= 
			x_{\alpha_1+\alpha_2}(N_{\alpha_1,\alpha_2}u).
		\end{align*} Thus, $u = \mu$. In the same way, regardless of whether $\alpha_3$ is long or short,
		\red{if we don't care for restricting $n$ to $\geq 3$, the outer commutator is redundant}
		\begin{align*}
			x_{\alpha_{\max}}(N_{\alpha_1+\alpha_2,\alpha_3,1}
				&N_{\alpha_1+\alpha_2+\alpha_3,\alpha_{\max}-\alpha_1-\alpha_2-\alpha_3} \nu)\\
			&=
			\phi(x_{\alpha_{\max}}(N_{\alpha_1+\alpha_2,\alpha_3,1}
				N_{\alpha_1+\alpha_2+\alpha_3,\alpha_{\max}-\alpha_1-\alpha_2-\alpha_3} \nu))\\
			&=
			[[\phi(x_{\alpha_1}(\mu) x_{\alpha_1+\alpha_2}(\nu)), \phi(x_{\alpha_3}(1))],
				\phi(x_{\alpha_{\max}-\alpha_1-\alpha_2-\alpha_3}(1))]\\
			&= 
			[[x_{\alpha_1}(u) x_{\alpha_1+\alpha_2}(v) x_{\alpha_{\max}}(*), x_{\alpha_3}(1)],
				x_{\alpha_{\max}-\alpha_1-\alpha_2-\alpha_3}(1)]\\
			&=
			x_{\alpha_{\max}}(N_{\alpha_1+\alpha_2,\alpha_3,1}	
				N_{\alpha_1+\alpha_2+\alpha_3,\alpha_{\max}-\alpha_1-\alpha_2-\alpha_3} v).
		\end{align*}				
		Therefore, $v = \nu$ and \eqref{eq:lemma:n:1:to:n:a1+a12} holds.
		
		Next, consider the image $y$ under $\phi$ of $x_{\alpha_{\max}-\alpha_1}(r) x_{\alpha_3+\dots+\alpha_n}(s)$.
		By Lemma \ref{lemma:simple:roots:as:C} and Corollary \ref{cor:centralizers}, the subgroup
		\[
			\begin{aligned}
				X_{\alpha_n} &\cdot
				X_{\alpha_2+\dots+\alpha_n}	\cdot
				X_{\alpha_2+\dots+\alpha_n+\dots+\alpha_3} \cdot
				X_{\alpha_3+\dots+\alpha_n}	\cdot
				X_{\alpha_2+\dots+\alpha_n+\dots+\alpha_2} \\
				&\cdot
				X_{\alpha_3+\dots+\alpha_n+\dots+\alpha_3} \cdot
				X_{\alpha_{\max}-\alpha_1} \cdot
				X_{\alpha_{\max}-\alpha_1-\alpha_2} \cdot
				X_{\alpha_1+\dots+\alpha_n} \cdot
				X_{\alpha_{\max}}.
			\end{aligned}
		\]
		is preserved by $\phi$. Together with the fact that $\phi$ preserves $\Up^{(n-2)}$, this yields
		\[
			\begin{aligned}
				y &= 
				x_{\alpha_2+\dots+\alpha_n}(p_{2n})	\cdot
				x_{\alpha_2+\dots+\alpha_n+\dots+\alpha_3}(p_{23}) \cdot
				x_{\alpha_3+\dots+\alpha_n}(p_{3n})	\\
				&\cdot
				x_{\alpha_2+\dots+\alpha_n+\dots+\alpha_2}(p_{22}) \cdot
				x_{\alpha_3+\dots+\alpha_n+\dots+\alpha_3}(p_{33}) \cdot
				x_{\alpha_{\max}-\alpha_1}(p_{12})\\
				&\cdot
				x_{\alpha_{\max}-\alpha_1-\alpha_2}(p_{13}) \cdot
				x_{\alpha_1+\dots+\alpha_n}(p_{1n}) \cdot
				x_{\alpha_{\max}}(*)
			\end{aligned}
		\]
		for some $p_{ij} \in F$. The clean up goes as usual. First, observe that
		\begin{equation}
			\label{eq:lemma:n:1:to:n:1}
			\begin{aligned}
				x_{\alpha_{\max}}(N_{\alpha_{\max}-\alpha_1, \alpha_1} r)
				&=
				\phi(x_{\alpha_{\max}}(N_{\alpha_{\max}-\alpha_1, \alpha_1} r))\\
				&=
				\phi([x_{\alpha_{\max}-\alpha_1}(r) x_{\alpha_3+\dots+\alpha_n}(s), x_{\alpha_1}(1)])\\
				&=
				[\phi(x_{\alpha_{\max}-\alpha_1}(r) x_{\alpha_3+\dots+\alpha_n}(s)),\phi(x_{\alpha_1}(1))]\\
				&= 
				[y, x_{\alpha_1}(1)] \\
				&= 
				x_{\alpha_1+\dots+\alpha_n}(\pm p_{2n})	\cdot
					x_{\alpha_1+\dots+\alpha_n+\dots+\alpha_3}(\pm p_{23}) \cdot
					x_{\alpha_{\max}-\alpha_1}(\pm p_{22}) \\
				&\quad\cdot
					x_{\alpha_{\max}}(N_{\alpha_{\max}-2\alpha_1,\alpha_1,2} p_{22}+
					N_{\alpha_{\max}-\alpha_1, \alpha_1} p_{12}).
			\end{aligned}
		\end{equation}
		Comparing the left-hand side of \eqref{eq:lemma:n:1:to:n:1} with the right-hand side thereof, we deduce that
		$p_{2n}=p_{23}=p_{22}=0$ and $p_{12}=r$.
		\red{In what follows we assume that $n \geq 4$:}
		Next,
		\begin{equation}
			\label{eq:lemma:n:1:to:n:2}
			\begin{aligned}
				x_{\alpha_1+\dots+\alpha_n}(N_{\alpha_2+\dots+\alpha_n, \alpha_1+\alpha_2} s)
				&=
				\phi(x_{\alpha_1+\dots+\alpha_n}(N_{\alpha_2+\dots+\alpha_n, \alpha_1+\alpha_2} s))\\
				&=
				\phi([x_{\alpha_{\max}-\alpha_1}(r) x_{\alpha_3+\dots+\alpha_n}(s), x_{\alpha_1+\alpha_2}(1)])\\
				&=
				[\phi(x_{\alpha_{\max}-\alpha_1}(r) x_{\alpha_3+\dots+\alpha_n}(s)),\phi(x_{\alpha_1+\alpha_2}(1))]\\
				&= 
				[y, x_{\alpha_1+\alpha_2}(1)] \\
				&= 
				x_{\alpha_1+\dots+\alpha_n}(N_{\alpha_2+\dots+\alpha_n, \alpha_1+\alpha_2} p_{3n})\\
				&\quad\cdot
				x_{\alpha_{\max}-\alpha_1-\alpha_2}(\pm p_{33}) \cdot
				x_{\alpha_{\max}}(\pm p_{33} \pm 2 p_{13})
			\end{aligned}
		\end{equation}
		Comparing the left-hand side and the right-hand side of \eqref{eq:lemma:n:1:to:n:2}, we conclude that
		$p_{33}=p_{13}=0$ and $p_{3n}=s$. Summing up,
		\begin{equation}
			\label{eq:lemma:n:1:to:n:a1,-2+a3n}
			\begin{aligned}
				\phi(x_{\alpha_{\max}-\alpha_1}(r) x_{\alpha_3+\dots+\alpha_n}(s)) &= 
				x_{\alpha_{\max}-\alpha_1}(r) x_{\alpha_3+\dots+\alpha_n}(s)
				x_{\alpha_1+\dots+\alpha_n}(p_{1n})
				x_{\alpha_{\max}}(*)
			\end{aligned}
		\end{equation}
		Now, combining \eqref{eq:lemma:n:1:to:n:a1+a12} and \eqref{eq:lemma:n:1:to:n:a1,-2+a3n} we have
		\begin{equation}
			\label{eq:lemma:n:1:to:n:construction}
			\begin{aligned}\relax
				\phi(x_{\alpha_1+\dots+\alpha_n}&(N_{\alpha_1+\alpha_2, \alpha_3+\dots+\alpha_n} \nu s)
				x_{\alpha_{\max}}(N_{\alpha_1, \alpha_{\max}-\alpha_1} \mu r))\\
				&=	
			    \phi([x_{\alpha_1}(\mu) x_{\alpha_1+\alpha_2}(\nu), 
			    	x_{\alpha_{\max}-\alpha_1}(r) x_{\alpha_3+\dots+\alpha_n}(s)])\\
				&=
			    [\phi(x_{\alpha_1}(\mu) x_{\alpha_1+\alpha_2}(\nu)), 
			    	\phi(x_{\alpha_{\max}-\alpha_1}(r) x_{\alpha_3+\dots+\alpha_n}(s))]\\
			    &=
			    [x_{\alpha_1}(\mu) x_{\alpha_1+\alpha_2}(\nu) x_{\alpha_{\max}}(*),
					x_{\alpha_{\max}-\alpha_1}(r) x_{\alpha_3+\dots+\alpha_n}(s)\\
				&\qquad\qquad\cdot
					x_{\alpha_1+\dots+\alpha_n}(p_{1n})
					x_{\alpha_{\max}}(*)]\\
				&=
				x_{\alpha_1+\dots+\alpha_n}(N_{\alpha_1+\alpha_2, \alpha_3+\dots+\alpha_n} \nu s)
				x_{\alpha_{\max}}(N_{\alpha_1, \alpha_{\max}-\alpha_1} \mu r)			
			\end{aligned}
		\end{equation}
		Set $\nu = (N_{\alpha_1+\alpha_2, \alpha_3+\dots+\alpha_n})^{-1}$ and
		$\mu = (N_{\alpha_1, \alpha_{\max}-\alpha_1} )^{-1}$. Then \eqref{eq:lemma:n:1:to:n:construction}
		rewrites as
		\begin{equation}
			\label{eq:lemma:n:1:to:n:construction'}
				\phi(x_{\alpha_1+\dots+\alpha_n}(s)x_{\alpha_{\max}}(r))
				=
				x_{\alpha_1+\dots+\alpha_n}(s)x_{\alpha_{\max}}(r)			
		\end{equation}
		for arbitrary $r,s \in F$. Finally, \eqref{eq:lemma:n:1:to:n:construction'} yields
		\begin{align*}
			x_{\alpha_1+\dots+\alpha_n}&(N_{\alpha_n,\alpha_1+\dots+\alpha_{n-1},1} \xi)
				x_{\alpha_{\max}}(N_{\alpha_n,\alpha_1+\dots+\alpha_{n-1},2} \xi+
				N_{\alpha_1+\dots+\alpha_n,\alpha_1+\dots+\alpha_{n-1}} \zeta)\\
			&=
			[\phi(x_{\alpha_n}(\xi) x_{\alpha_1+\dots+\alpha_n}(\zeta)), \phi(x_{\alpha_1+\dots+\alpha_{n-1}}(1))]\\
			&=
			[x_{\alpha_n}(\xi) x_{\alpha_1+\dots+\alpha_n}(b) x_{\alpha_{\max}}(*), x_{\alpha_1+\dots+\alpha_{n-1}}(1)]\\					&= 
			x_{\alpha_1+\dots+\alpha_n}(N_{\alpha_n,\alpha_1+\dots+\alpha_{n-1},1} \xi)\\
			&\qquad\qquad\cdot
				x_{\alpha_{\max}}(N_{\alpha_n,\alpha_1+\dots+\alpha_{n-1},2} \xi+
				N_{\alpha_1+\dots+\alpha_n,\alpha_1+\dots+\alpha_{n-1}} b).
		\end{align*}
		Thus $N_{\alpha_1+\dots+\alpha_n,\alpha_1+\dots+\alpha_{n-1}} \zeta = N_{\alpha_1+\dots+\alpha_n,\alpha_1+\dots+\alpha_{n-1}} b$ and as $2F=F$, $\zeta=b$. This completes the proof.
	\end{proof}
\end{lemma}

\begin{lemma}
	\label{lemma:field:autom}
	Assume $\rk(\Phi) \geq 3$ and $2F=F$.
	Set $\Pi' = \{ \alpha_1, \dots, \alpha_{n-1}, \alpha_{n-1}+\alpha_n \}$.
	Let $\phi \in \PC(\Up)$ satisfy the following conditions:
	\begin{enumerate}
		\item $\phi(X_{\alpha_n}) \in X_n \cdot X_{\alpha_1+\dots+\alpha_n} \cdot X_{\alpha_{\max}}$;
		\item $\phi(x_{\alpha}(1)) = x_\alpha(1) \cdot x_{\alpha_{\max}}(*)$ for all $\alpha \in \Pi'$;
		\item $\phi(x_{\alpha}(1)) = x_\alpha(1)$ for all $\alpha \in \Phi^S \setminus \Pi'$;
		\item $\phi(x_{\alpha}(2)) = x_\alpha(2)$ for all $\alpha \in \Phi^L \setminus \{\alpha_n\}$.
	\end{enumerate}
	Then there exists a field automorphism $\tau$ and a central $\PC$-map $C$ such that
	\[
		\tau \circ C \circ \phi (x_\alpha(\xi)) = x_\alpha(\xi),
	\]
	for all $\alpha \neq \alpha_{n-1}+\alpha_n$ and
	\[
		\tau \circ C \circ \phi (x_{\alpha_{n-1}+\alpha_n}(\xi)) = x_{\alpha_{n-1}+\alpha_n}(\xi) x_{\alpha_{\max}}(*).
	\]
	\begin{proof}
		Let $f_i(\xi) : F \rightarrow F$ be a function such that
		\[
			\phi(x_{\alpha_i}(\xi)) = x_{\alpha_i}(*) x_{\alpha_{\max}}(-f_i(\xi))
		\]
		for $1 \leq i \leq n-1$ and 
		\[
			\phi(x_{\alpha_n}(\xi)) = x_{\alpha_n}(*) x_{\alpha_1+\dots+\alpha_n}(*) x_{\alpha_{\max}}(-f_n(\xi)).
		\]
		Let $C$ be the standard central $\PC$-map defined by the map
		$f(t_1,\dots,t_n)=f_1(t_1)+\dots+f_n(t_n)$ (cf. Lemma \ref{lemma:standard:central:map}). Clearly,
		\[
			C \circ \phi (x_{\alpha}) \in X_\alpha
		\]
		for all $\alpha \in \Pi \setminus \{\alpha_n\}$ and $C \circ \phi (x_{\alpha_n}) \in X_{\alpha_n}
		\cdot X_{\alpha_1+\dots+\alpha_n}$. Replace $\phi$ with $C \circ \phi$.
		 
		Now define $\tau : F \rightarrow F$ so that
		\[
			\phi(x_{\alpha_1}(\xi)) = x_{\alpha_1}(\tau(\xi))
		\]
		for all $\xi \in F$. Observe that 
		\[
			\phi(x_{\alpha_1+\alpha_2}(\xi)) = \phi([x_{\alpha_1}(\xi),x_{\alpha_2}(1)])
			= [x_{\alpha_1}(\tau(\xi)), x_{\alpha_2}(1)] = x_{\alpha_1+\alpha_2}(\tau(\xi))
		\]	
		On the other hand, 
		\[
			x_{\alpha_1+\alpha_2}(\zeta) = [x_{\alpha_1}(1), x_{\alpha_2}(\zeta)] = \phi(x_{\alpha_1+\alpha_2}(\xi)) = x_{\alpha_1+\alpha_2}(\tau(\xi)),
		\]
		where $\phi(x_{\alpha_2}(\xi))=\zeta$, thus $\zeta = \tau(\xi)$ and $\phi(x_{\alpha_2}(\xi))=x_{\alpha_2}(\tau(\xi))$.
		In the same way, \red{HAVE TO WRITE THIS EXPLICITLY}
		\[
			\phi(x_{\alpha}(\xi)) = x_{\alpha}(\tau(\xi))
		\]
		for all $\alpha \neq \alpha_n$ and
		\[
			\phi(x_{\alpha_{n-1}+\alpha_n}(\xi)) = x_{\alpha_{n-1}+\alpha_n}(\tau(\xi)) x_{\alpha_{\max}}(*).
		\]
		
		Next, we will show that $\tau$ is a field automorphism. Indeed, multiplicativity follows from
		\begin{align*}
			x_{\alpha_1+\alpha_2}(\tau(\xi\zeta)) &= \phi(x_{\alpha_1+\alpha_2}(\xi\zeta)) 
			= \phi([x_{\alpha_1}(\xi),x_{\alpha_2}(\zeta)]) \\
			&=
			[x_{\alpha_1}(\tau(\xi)),x_{\alpha_2}(\tau(\zeta))] = x_{\alpha_1+\alpha_2}(\tau(\xi)\tau(\zeta)).
		\end{align*}
		\red{Mistake in Chen's paper here}
		the proof of additivity of $\tau$ requires a more sophisticated analysis. First, we will partially compute
		the image $y$ of $x_{\alpha_1}(a)x_{\alpha_2}(b)$ under $\phi$. 
		Note that $[y,x_{\alpha_{\max}-\alpha_1}(1)]=x_{\alpha_{\max}}(2 y_{12})$. Thus,
		\begin{align*}
			x_{\alpha_{\max}}(\tau(2)\tau(a)) = \phi(x_{\alpha_{\max}}(2 a)) &= 
			\phi([x_{\alpha_1}(a)x_{\alpha_2}(b),x_{\alpha_{\max}-\alpha_1}(1)])\\
			&= 
			[y, x_{\alpha_{\max}-\alpha_1}(1)] = x_{\alpha_{\max}}(2 y_{12}).
		\end{align*}
		Consequently, $y_{12} = \frac{\tau(2)}{2} \tau(a)$. In the same way, an explicit computation shows that
		\[
			[y,x_{\alpha_{\max}-\alpha_1-\alpha_2}(1)] = 
			x_{\alpha_{\max}}(2 y_{13}) x_{\alpha_{\max}-\alpha_2}(y_{23}).
		\]		
		Therefore,
		\begin{align*}
		  x_{\alpha_{\max}-\alpha_1}(\tau(b))
		  =	
	      \phi(x_{\alpha_{\max}-\alpha_1}(b))
		  &=
		  [\phi(x_{\alpha_1}(a)x_{\alpha_2}(b)),\phi(x_{\alpha_{\max}-\alpha_1-\alpha_2}(1))]\\
		  &=
	      [y,x_{\alpha_{\max}-\alpha_1-\alpha_2}(1)]=
		  x_{\alpha_{\max}}(2 y_{13}) x_{\alpha_{\max}-\alpha_2}(y_{23}).
		\end{align*}
		Thus $y_{13}=0$ and $y_{23}=\tau(b)$.
		Finally, by substituting $b=0$ we get $y = x_{\alpha_1}(\tau(a))$ and thus
		$\tau(a) = \frac{\tau(2)}{2} \tau(a)$, thus $\tau(2)=2$ and $y_{12}=\tau(a)$.
		
		Now consider
		\[
			x_{\alpha_1+\alpha_2}(\tau(b-a))=\phi(x_{\alpha_1+\alpha_2}(b-a)) = 
			\phi([x_{\alpha_1}(1)x_{\alpha_2}(a),x_{\alpha_1}(1)x_{\alpha_2}(b)])
			= [y,z],
		\]
		with $y_{12}=1,y_{23}=\tau(a),z_{12}=1,z_{23}=\tau(b)$. Therefore
		$[y,z]_{13}=\tau(a)-\tau(b)$. Thus 
		\[
			\tau(a)-\tau(b)=x_{\alpha_1+\alpha_2}(\tau(b-a))_{13} = \tau(b-a).
		\]
		Let $\psi$ be the field automorphism of $\Up$ corresponding to $\tau^{-1}$. Set $\phi'$ with
		$\psi \circ C \circ \phi$. We have shown that
		\begin{equation}
			\label{eq:lemma:field:autom:final}
			\begin{aligned}
				\phi'(x_{\alpha_n}(\xi)) &\in X_{\alpha_n} \cdot X_{\alpha_1+\dots+\alpha_n} \\
				\phi'(x_{\alpha_{n-1}+\alpha_n}(\xi)) &= x_{\alpha_{n-1}+\alpha_n}(\xi) \cdot x_{\alpha_{\max}}(*) \\			
				\phi'(x_{\alpha}(\xi)) &= x_{\alpha}(\xi)
			\end{aligned}
		\end{equation}
		for all $\alpha \neq \alpha_n, \alpha_{n-1}+\alpha_n$.				

		Finally, by Lemma \ref{lemma:n:1:to:n} we get 
		\[
			\phi'(x_{\alpha_n}(\xi)) = x_{\alpha_n}(\xi) x_{\alpha_{\max}}(t).
		\]
		On the other hand, \eqref{eq:lemma:field:autom:final} requires that $t=0$. 
	\end{proof}
\end{lemma}

We have just proved that any $\PC$ map, up to standard ones is an almost identity map. The root
$\alpha_n+\alpha_{n-1}$ is exceptional for despite being of height 2 is not a commutator.

\begin{theorem}
	\label{theorem:up:to:ai}
	Assume $n \geq 4$ and $6F=F$. Then every $\phi \in \PC(\Up)$ can be presented in the following form:
	\[
		\phi = \Int_C \circ E \circ Int_{C'} \circ D \circ Q \circ Z \circ \tau \circ A,
	\]
	where $\Int_C$ and $\Int_{C'}$ are inner, $E$ an extremal, $D$ a diagonal, $Q$ a quasi-diagonal,
	$\tau$ a field automorphisms; $Z$ -- a central $\PC$-map and $A$ -- an almost identity $\PC$-map.
	\begin{proof}
		Sequentially apply Lemmas \ref{lemma:T12}, \ref{lemma:U1}, \ref{lemma:simple:roots}, 
		\ref{lemma:X:dirty} and \ref{lemma:field:autom}.
	\end{proof}
\end{theorem}

\section{Almost identity is central}

In this section we will characterize almost identity maps.

\begin{lemma}
	\label{lemma:ai:giblets}
	Assume $\rk(\Phi) \geq 2$ and $2F=F$.
	Let $\phi$ be an almost identity $\PC$-map on $\Up$ and $a \in \Up$. 
	Then $\phi(a)_{ij} = a_{ij}$ for all $2 \leq i < j \le 2n-1$.
	\begin{proof}
		Fix $i,j$ such that $2 \leq i < j \le 2n-1$.
		By Lemma \ref{lemma:extract:middle} there exist roots $\beta,\gamma \in \Phi$ such that
		\[
			x_{\alpha_{\max}}(\pm 2 a_{ij}) = [x_\beta(1),[x_\gamma(1),a]].
		\]
		Then, as $\phi$ is an almost identity map,
		\begin{align*}
			x_{\alpha_{\max}}(\pm 2 a_{ij})
			&=
			\phi(x_{\alpha_{\max}}(\pm 2 a_{ij}))
			=
			\phi([x_\beta(1),[x_\gamma(1),a]])\\
			&=
			[\phi(x_\beta(1)),[\phi(x_\gamma(1)),\phi(a)]]\\
			&=
			[x_\beta(1) x_{\alpha_{\max}}(*),[x_\gamma(1) x_{\alpha_{\max}}(*),\phi(a)]]\\
			&=
			[x_\beta(1),[x_\gamma(1), \phi(a)]]
			=
			x_{\alpha_{\max}}(\pm 2 \phi(a)_{ij}).	
		\end{align*}		
		Thus, $\phi(a)_{ij} = a_{ij}$.
	\end{proof}
\end{lemma}

\begin{lemma}
	\label{lemma:ai:only:skin}
	Let $\phi$ be an almost identity $\PC$-map on $\Up$ and $a \in \U_1$, i.e 
	\begin{equation}
		\label{eq:lemma:ai:only:skin:statement}
		a = \prod_{\alpha \in \Phi, m_1(\alpha)=1} x_{\alpha}(a_\alpha) \cdot x_{\alpha_{\max}}(a_{\alpha_{\max}}).
	\end{equation}
	Then $\phi(a)_{ij} = a_{ij}$ whenever $(i,j) \neq (1,2n)$.
	\begin{proof}
		We have already shown that $\phi(\U_1) \subseteq \U_1$ for any $\PC$-map. Thus
		\[
			\phi(a) = \prod_{\alpha \in \Phi, m_1(\alpha)=1} 
				x_{\alpha}(\phi(a)_\alpha) \cdot x_{\alpha_{\max}}(a_{\alpha_{\max}}),
		\] 
		where the product is taken in the same order as in \eqref{eq:lemma:ai:only:skin:statement}.
		It is enough to show that $a_\alpha = \phi(a)_\alpha$ for all $\alpha \in \Phi$
		such that $m_1(\alpha) = 1$. Observe that for any short root $\alpha \in \Phi$ such 
		that $m_1(\alpha) = 1$ there
		exists the unique short root $\beta$ such that $m_1(\beta)=1$ and $\alpha+\beta \in \Phi$.
		Specifically, $\beta = \alpha_{\max}-\alpha$. Thus for any $\alpha \in \Phi$ 
		such that $m_1(\alpha)=1$ we have
		\[
			x_{\alpha_{\max}}(\pm 2 a_\alpha) = \phi(x_{\alpha_{\max}}(\pm 2 a_\alpha)) 
			= \phi([a, x_{\alpha_{\max}-\alpha}(1)])
			= \phi(x_{\alpha_{\max}}(\pm 2 \phi(a)_\alpha).
		\]
		Thus $\phi(a)_{\alpha} = a_{\alpha}$ for all $\alpha$ such that $m_1(\alpha)=1$.
	\end{proof}
\end{lemma}

\begin{lemma}
	\label{lemma:ai:prod:of:simple}
	Assume $\rk \Phi \geq 4$ and $2F=F$.
	Let $\phi$ be an almost identity $\PC$-map on $\Up$ and 
	\[
		a = \prod_{i=1}^n x_{\alpha_i}(\xi_i) \in \Up.
	\]
	Then $\phi(a) = a x_{\alpha_{\max}}(*)$.
	\begin{proof}
		By Lemma \ref{lemma:ai:giblets} we have $a_{ij} = \phi(a)_{ij}$ for all
		$2 \leq i < j \leq 2n-1$. It is easily seen that
		\[
			\phi(a) = \prod_{i=2}^n x_{\alpha_i}(\xi_i) \cdot
			\prod_{\alpha \in \Phi, m_1(\alpha) \geq 1} x_\alpha(\zeta_\alpha).
		\]
		for some $\zeta_\alpha \in F$. 
		It is enough for us to show that $\zeta_{\alpha_1}=\xi_1$ and $\zeta_{\alpha}=0$ for
		all $\alpha \in \Phi$, $m_1(\alpha) = 1$, $\alpha \neq \alpha_1$.
		Pick a root $\alpha$ be such that $m_1(\alpha)=1$. We continue case by case.
		
		Let $1 \leq j \leq n-2$. Consider
		\begin{equation}
			\begin{aligned}
				x_{\alpha_1+\dots+\alpha_j+\alpha_{j+1}}&(N_{\alpha_{j+1},\alpha_1+\dots+\alpha_j} \xi_{{j+1}}) \cdot
				x_{\alpha_{\max}}(\pm 2 \zeta_{\alpha_{\max}-\alpha_1-\dots-\alpha_j}) \\
				&=
				[\phi(a), x_{\alpha_1+\dots+\alpha_j}(1)]
				=
				\phi([a, x_{\alpha_1+\dots+\alpha_j}(1)]) \\
				&= 
				\phi(x_{\alpha_1+\dots+\alpha_j+\alpha_{j+1}}(N_{\alpha_{j+1},\alpha_1+\dots+\alpha_j} \xi_{{j+1}}))\\
				&= 
				x_{\alpha_1+\dots+\alpha_j+\alpha_{j+1}}(N_{\alpha_{j+1},\alpha_1+\dots+\alpha_j} \xi_{{j+1}}).
			\end{aligned}
		\end{equation}
		Thus $\zeta_{\alpha_{\max}-\alpha_1-\dots-\alpha_j} = 0$ for all $1 \leq j \leq n-2$.
		
		Next, let $2 \leq j \leq n-1$ and observe that
		\begin{equation}
			\begin{aligned}
				x_{\alpha_{\max}-\alpha_1-\dots-\alpha_{j-1}}&
					(N_{\alpha_{j},\alpha_{\max}-\alpha_1-\dots-\alpha_j}\xi_j) \cdot
				x_{\alpha_{\max}}(\pm 2 \zeta_{\alpha_1+\dots+\alpha_j})\\
				&=
				[\phi(a), x_{\alpha_{\max}-\alpha_1-\dots-\alpha_j}(1)]
				=
				\phi([a, x_{\alpha_{\max}-\alpha_1-\dots-\alpha_j}(1)]) \\
				&= 
				\phi(x_{\alpha_{\max}-\alpha_1-\dots-\alpha_{j-1}}
					(N_{\alpha_{j},\alpha_{\max}-\alpha_1-\dots-\alpha_j} \xi_j))\\
				&= x_{\alpha_{\max}-\alpha_1-\dots-\alpha_{j-1}}
					(N_{\alpha_{j},\alpha_{\max}-\alpha_1-\dots-\alpha_j}\xi_j).
			\end{aligned}
		\end{equation}
		Thus $\zeta_{\alpha_1+\dots+\alpha_j} = 0$ for all $2 \leq j \leq n-1$.
		
		It is only deft to deal with $\zeta_{\alpha_1}$ and $\zeta_{\alpha_1+\dots+\alpha_n}$.
		The former is handles in the following way:
		\red{This step heavily utilizes $n \geq 4$}
		\begin{equation}
			\begin{aligned}
				x_{\alpha_1+\alpha_2+\alpha_3}(N_{\alpha_1+\alpha_2,\alpha_3} N_{\alpha_1,\alpha_2} \zeta_{\alpha_1})
				&=
				[x_{\alpha_1+\alpha_2}(N_{\alpha_1,\alpha_2} \zeta_{\alpha_1})
				x_{\alpha_2+\alpha_3}(N_{\alpha_3,\alpha_2} \xi_3), x_{\alpha_3}(1)]\\
				&=
				[[\phi(a), x_{\alpha_2}(1)],x_{\alpha_3}(1)]
				=
				\phi([[a,x_{\alpha_2}(1)],x_{\alpha_3}(1)])\\
				&= 
				\phi([x_{\alpha_1+\alpha_2}(N_{\alpha_1,\alpha_2} \xi_1)
				x_{\alpha_2+\alpha_3}(N_{\alpha_3,\alpha_2}\xi_3),x_{\alpha_3}(1)])\\
				&=
				\phi(x_{\alpha_1+\alpha_2+\alpha_3}(N_{\alpha_1+\alpha_2,\alpha_3} N_{\alpha_1,\alpha_2} \xi_1))\\
				&=
				x_{\alpha_1+\alpha_2+\alpha_3}(N_{\alpha_1+\alpha_2,\alpha_3} N_{\alpha_1,\alpha_2} \xi_1)
			\end{aligned}
		\end{equation}
		Thus $\zeta_{\alpha_1}=\xi_1$. Finally, we have to show that 
		$\zeta_{\alpha_1+\dots+\alpha_n} = 0$. This case requires heavy artillery, Lemma \ref{lemma:n:1:to:n}.
		Indeed,
		\begin{align*}
			x_{\alpha_1+\dots+\alpha_n}&(N_{\alpha_n,\alpha_1+\dots+\alpha_{n-1},1} \xi_n)
				x_{\alpha_{\max}}(N_{\alpha_n,\alpha_1+\dots+\alpha_{n-1},2} \xi_n)\\
			&= 
			\phi(x_{\alpha_1+\dots+\alpha_n}(N_{\alpha_n,\alpha_1+\dots+\alpha_{n-1},1} \xi_n)
				x_{\alpha_{\max}}(N_{\alpha_n,\alpha_1+\dots+\alpha_{n-1},2} \xi_n))\\
			&= 
			\phi([a, x_{\alpha_1+\dots+\alpha_{n-1}}(1)])
			[\phi(a), \phi(x_{\alpha_1+\dots+\alpha_{n-1}}(1))]\\
			&=
			[\prod_{i=1}^n x_{\alpha_i}(\xi_i) \cdot
			x_{\alpha_1+\dots+\alpha_n}(\zeta_{\alpha_1+\dots+\alpha_n}) \cdot
			x_{\alpha_{\max}}(*), \phi(x_{\alpha_1+\dots+\alpha_{n-1}}(1))]\\
			&=
			x_{\alpha_1+\dots+\alpha_n}(N_{\alpha_n,\alpha_1+\dots+\alpha_{n-1},1} \xi_n)\\
			&\quad\cdot
			x_{\alpha_{\max}}(N_{\alpha_n,\alpha_1+\dots+\alpha_{n-1},2} \xi_n + 
				N_{\alpha_1+\dots+\alpha_n,\alpha_1+\dots+\alpha_{n-1}} \zeta_{\alpha_1+\dots+\alpha_n}),
		\end{align*}
		where the first equality is by Lemma \ref{lemma:n:1:to:n}.
		Thus $\zeta_{\alpha_1+\dots+\alpha_n} = 0$.
	\end{proof}
\end{lemma}

\begin{lemma}
	\label{lemma:ai:skin:max}
	Assume $n \geq 4$ and $2F=F$.
	Let $\phi$ be an almost identity $\PC$-map on $\Up$ and $a \in \U_1^{(2)}$. Then $\phi(a)=a$.
	\begin{proof}
		As, $a \in \U^{(2)}$ we can decompose $a$ as follows:
		\begin{equation}
			\label{eq:lemma:ai:skin:max:1}
			a = \left(\prod_{i=2}^{n-1}
					x_{\alpha_1+\dots+\alpha_i}(u_i) x_{\alpha_{\max}-\alpha_1-\dots-\alpha_{i-1}}(v_i)\right) 
				    x_{\alpha_1+\dots+\alpha_n}(u_n) 
					x_{\alpha_{\max}}(v_1).
		\end{equation}
		In order to show that $a$ is preserved by $\phi$, 
		we will present $a$ as a commutator of matrices already known to be preserved by $\phi$. 
		Let 
		\begin{align*}
			b &= \prod_{i=1}^{n-1} x_{\alpha_1+\dots+\alpha_i}(\zeta_i) x_{\alpha_{\max}-\alpha_1-\dots-\alpha_i}(\eta_i)
			&
			c &= \prod_{i=1}^n x_{\alpha_i}(\xi_i).
		\end{align*}				
		A direct calculation shows that
		\begin{equation}
			\label{eq:lemma:ai:skin:max:big}
			\begin{aligned}\relax
				[b, c] 
				&= 
				\prod_{i=1}^n {}^{\prod_{j=1}^{i-1} x_{\alpha_j}(\xi_j)} [b, x_{\alpha_i}(\xi_i)]\\
				&= 
				\prod_{i=1}^n  \left(\left[\prod_{j=1}^{i-1} x_{\alpha_j}(\xi_j), [b, x_{\alpha_i}(\xi_i)]\right]
					[b, x_{\alpha_i}(\xi_i)]\right).
			\end{aligned}
		\end{equation}
		We will compute \eqref{eq:lemma:ai:skin:max:big} step by step. First, observe that
		\begin{equation}
			\label{eq:lemma:ai:skin:max:b:xai}
			\begin{aligned}\relax
				[b, x_{\alpha_i}(\xi_i)] = 
				\begin{cases}
					x_{\alpha_{\max}}(N_{\alpha_{\max}-\alpha_1,\alpha_1} \eta_1 \xi_1), & i=1; \\
					\begin{aligned}
						&x_{\alpha_1+\dots+\alpha_i}(N_{\alpha_1+\dots+\alpha_{i-1},\alpha_i} \zeta_{i-1} \xi_i)\\
						&\quad\cdot x_{\alpha_{\max}-\alpha_1-\dots-\alpha_{i-1}}
							(N_{\alpha_{\max}-\alpha_1-\dots-\alpha_i,\alpha_i} \eta_i \xi_i)
					\end{aligned}, & 1<i<n; \\
					\begin{aligned}
						&x_{\alpha_1+\dots+\alpha_n}(N_{\alpha_1+\dots+\alpha_{n-1},\alpha_n,1} \zeta_{n-1} \xi_n)\\
						&\quad\cdot x_{\alpha_{\max}}(N_{\alpha_1+\dots+\alpha_{n-1},\alpha_n,2} \zeta_{n-1}^2 \xi_n)
					\end{aligned}, & i=n.
				\end{cases}
			\end{aligned}
		\end{equation}
		Next, in a view of \eqref{eq:lemma:ai:skin:max:b:xai}, if $1 < i < n$ we have
		\begin{equation}
			\label{eq:lemma:ai:skin:max:prodaj:b:xai}
			\begin{aligned}\relax
				[\prod_{j=1}^{i-1} x_{\alpha_j}(\xi_j), &[b, x_{\alpha_i}(\xi_i)]] \\
				&= 
				[\prod_{j=1}^{i-1} x_{\alpha_j}(\xi_j), 
					x_{\alpha_1+\dots+\alpha_i}(N_{\alpha_1+\dots+\alpha_{i-1},\alpha_i} \zeta_{i-1} \xi_i) \\
				&\qquad\qquad\cdot
					x_{\alpha_{\max}-\alpha_1-\dots-\alpha_{i-1}}
					(N_{\alpha_{\max}-\alpha_1-\dots-\alpha_i,\alpha_i} \eta_i \xi_i)] \\
				&= 
				[\prod_{j=1}^{i-1} x_{\alpha_j}(\xi_j), 
					x_{\alpha_1+\dots+\alpha_i}(N_{\alpha_1+\dots+\alpha_{i-1},\alpha_i} \zeta_{i-1} \xi_i)]\\
				&\quad\cdot
				{}^{x_{\alpha_1+\dots+\alpha_i}(*)} [\prod_{j=1}^{i-1} x_{\alpha_j}(\xi_j), 
					x_{\alpha_{\max}-\alpha_1-\dots-\alpha_{i-1}}	
						(N_{\alpha_{\max}-\alpha_1-\dots-\alpha_i,\alpha_i} \eta_i \xi_i)] \\
				&= 
				[\prod_{j=1}^{i-1} x_{\alpha_j}(\xi_j), 
					x_{\alpha_{\max}-\alpha_1-\dots-\alpha_{i-1}}	
						(N_{\alpha_{\max}-\alpha_1-\dots-\alpha_i,\alpha_i} \eta_i \xi_i)].
			\end{aligned}
		\end{equation}
		In the same way for $i=n$ we have
		\begin{equation}
			\label{eq:lemma:ai:skin:max:prodaj:b:xan}
			\begin{aligned}\relax
				[\prod_{j=1}^{n-1} x_{\alpha_j}(\xi_j), [b, x_{\alpha_n}(\xi_n)]] 
				&= 
				[\prod_{j=1}^{n-1} x_{\alpha_j}(\xi_j), 
					x_{\alpha_1+\dots+\alpha_n}(N_{\alpha_1+\dots+\alpha_{n-1},\alpha_n,1} \zeta_{n-1} \xi_n)\\
				&\qquad\qquad\cdot
				x_{\alpha_{\max}}(N_{\alpha_1+\dots+\alpha_{n-1},\alpha_n,2} \zeta_{n-1}^2 \xi_n)]\\					
				&= 
				[\prod_{j=1}^{n-1} x_{\alpha_j}(\xi_j), 
					x_{\alpha_1+\dots+\alpha_n}(N_{\alpha_1+\dots+\alpha_{n-1},\alpha_n,1} \zeta_{n-1} \xi_n)].
			\end{aligned}
		\end{equation}
		The next step is computing 
		\[
			A_{i}(t) = [\prod_{j=1}^{i} x_{\alpha_j}(\xi_j), x_{\alpha_{\max}-\alpha_1-\dots-\alpha_{i}}(t)].
		\] for $1 \leq i \leq n-1$. Observe that
		\begin{equation}
			\label{eq:lemma:ai:skin:max:Ai}
			\begin{aligned}\relax
				A_{i}(t) &= [\prod_{j=1}^{i} x_{\alpha_j}(\xi_j), x_{\alpha_{\max}-\alpha_1-\dots-\alpha_{i}}(t)]\\
				&= {}^{\prod_{j=1}^{i-1} x_{\alpha_j}(\xi_j)}[x_{\alpha_i}(\xi_i), 
					x_{\alpha_{\max}-\alpha_1-\dots-\alpha_{i}}(t)] 
				\cdot [\prod_{j=1}^{i-1} x_{\alpha_j}(\xi_j), x_{\alpha_{\max}-\alpha_1-\dots-\alpha_{i}}(t)] \\
				&= {}^{\prod_{j=1}^{i-1} x_{\alpha_j}(\xi_j)}[x_{\alpha_i}(\xi_i), 
					x_{\alpha_{\max}-\alpha_1-\dots-\alpha_{i}}(t)]\\
				&= 
				{}^{\prod_{j=1}^{i-1} x_{\alpha_j}(\xi_j)}
				x_{\alpha_{\max}-\alpha_1-\dots-\alpha_{i-1}}(N_{\alpha_i,\alpha_{\max}-\alpha_1-\dots-\alpha_i}\xi_i t)\\
				&=
				[\prod_{j=1}^{i-1} x_{\alpha_j}(\xi_j),
				x_{\alpha_{\max}-\alpha_1-\dots-\alpha_{i-1}}(N_{\alpha_i,\alpha_{\max}-\alpha_1-\dots-\alpha_i}\xi_i t)]\\
				&\quad\cdot
				x_{\alpha_{\max}-\alpha_1-\dots-\alpha_{i-1}}(N_{\alpha_i,\alpha_{\max}-\alpha_1-\dots-\alpha_i}\xi_i t)\\
				&=
				A_{i-1}(N_{\alpha_i,\alpha_{\max}-\alpha_1-\dots-\alpha_i}\xi_i t)
				\cdot
				x_{\alpha_{\max}-\alpha_1-\dots-\alpha_{i-1}}(N_{\alpha_i,\alpha_{\max}-\alpha_1-\dots-\alpha_i}\xi_i t).
			\end{aligned}
		\end{equation}
		Applying \eqref{eq:lemma:ai:skin:max:Ai} repeatedly to itself, we get
		\begin{equation}
			\label{eq:lemma:ai:skin:max:Ai'}
			\begin{aligned}\relax
				A_{i}(t) &= 
				A_{i-1}(N_{\alpha_i,\alpha_{\max}-\alpha_1-\dots-\alpha_i}\xi_i t)
				\cdot
				x_{\alpha_{\max}-\alpha_1-\dots-\alpha_{i-1}}(N_{\alpha_i,\alpha_{\max}-\alpha_1-\dots-\alpha_i}\xi_i t)\\
				&=
				A_{i-2}(N_{\alpha_{i-1},\alpha_{\max}-\alpha_1-\dots-\alpha_{i-1}} \xi_{i-1} 
					N_{\alpha_i,\alpha_{\max}-\alpha_1-\dots-\alpha_i}\xi_i t)\\
				&\quad\cdot					
				x_{\alpha_{\max}-\alpha_1-\dots-\alpha_{i-2}}
					(N_{\alpha_{i-1},\alpha_{\max}-\alpha_1-\dots-\alpha_{i-1}} \xi_{i-1} 
						N_{\alpha_i,\alpha_{\max}-\alpha_1-\dots-\alpha_i}\xi_i t)\\
					&\quad\cdot
				x_{\alpha_{\max}-\alpha_1-\dots-\alpha_{i-1}}(N_{\alpha_i,\alpha_{\max}-\alpha_1-\dots-\alpha_i}\xi_i t)\\
				&=\dots\\
				&= \prod_{j=1}^{i} x_{\alpha_{\max}-\alpha_1-\dots-\alpha_{j-1}}\left(\left(\prod_{k=j}^i N_{\alpha_k,\alpha_{\max}-\alpha_1-\dots-\alpha_k} \xi_k\right) \cdot t \right).
			\end{aligned}
		\end{equation}
		Substituting \eqref{eq:lemma:ai:skin:max:b:xai}, \eqref{eq:lemma:ai:skin:max:prodaj:b:xai} and \eqref{eq:lemma:ai:skin:max:prodaj:b:xan}
		into \eqref{eq:lemma:ai:skin:max:big} we get
		\begin{equation}
			\label{eq:lemma:ai:skin:max:big'}
			\begin{aligned}\relax
				[b,c] 
				&= 
				[b, x_{\alpha_1}(\xi_1)] \cdot
				\prod_{i=1}^n  \left(\left[\prod_{j=1}^{i-1} x_{\alpha_j}(\xi_j), [b, x_{\alpha_i}(\xi_i)]\right]
					\cdot [b, x_{\alpha_i}(\xi_i)]\right) \\
				&\quad\cdot
				\left(\left[\prod_{j=1}^{n-1} x_{\alpha_j}(\xi_j), [b, x_{\alpha_n}(\xi_n)]\right]
					\cdot [b, x_{\alpha_n}(\xi_n)]\right)\\
				&= 
				x_{\alpha_{\max}}(N_{\alpha_{\max}-\alpha_1,\alpha_1} \eta_1 \xi_1) \cdot
				\prod_{i=2}^{n-1}( A_{i-1}(N_{\alpha_{\max}-\alpha_1-\dots-\alpha_i,\alpha_i} \eta_i \xi_i)\\
				&\qquad\qquad\cdot
					x_{\alpha_1+\dots+\alpha_i}(N_{\alpha_1+\dots+\alpha_{i-1},\alpha_i} \zeta_{i-1} \xi_i)\\
				&\qquad\qquad\cdot
					x_{\alpha_{\max}-\alpha_1-\dots-\alpha_{i-1}}
						(N_{\alpha_{\max}-\alpha_1-\dots-\alpha_i,\alpha_i} \eta_i \xi_i))\\
				&\quad\cdot
				(A_{n-1}(N_{\alpha_1+\dots+\alpha_{n-1},\alpha_n,1} \zeta_{n-1} \xi_n)
					\cdot x_{\alpha_1+\dots+\alpha_n}(N_{\alpha_1+\dots+\alpha_{n-1},\alpha_n,1} \zeta_{n-1} \xi_n)\\
				&\qquad\qquad\cdot
					x_{\alpha_{\max}}(N_{\alpha_1+\dots+\alpha_{n-1},\alpha_n,2} \zeta_{n-1}^2 \xi_n)).
			\end{aligned}
		\end{equation}
		Next, substituting \eqref{eq:lemma:ai:skin:max:Ai'} into \eqref{eq:lemma:ai:skin:max:big'} we get
		\begin{equation}
			\label{eq:lemma:ai:skin:max:big''}
			\begin{aligned}\relax
				[b,c] 
				&= 
				x_{\alpha_{\max}}(N_{\alpha_{\max}-\alpha_1,\alpha_1} \eta_1 \xi_1)\\
				&\quad \cdot
				\prod_{i=2}^{n-1}( \prod_{j=1}^{i-1} x_{\alpha_{\max}-\alpha_1-\dots-\alpha_{j-1}}\left(-\eta_i \prod_{k=j}^i N_{\alpha_k,\alpha_{\max}-\alpha_1-\dots-\alpha_k} \xi_k\right)\\
				&\qquad\qquad\cdot
					x_{\alpha_1+\dots+\alpha_i}(N_{\alpha_1+\dots+\alpha_{i-1},\alpha_i} \zeta_{i-1} \xi_i)\\
				&\qquad\qquad\cdot
					x_{\alpha_{\max}-\alpha_1-\dots-\alpha_{i-1}}
						(N_{\alpha_{\max}-\alpha_1-\dots-\alpha_i,\alpha_i} \eta_i \xi_i))\\
				&\quad\cdot
				(\prod_{j=1}^{n-1} x_{\alpha_{\max}-\alpha_1-\dots-\alpha_{j-1}}\left(- \zeta_{n-1} \cdot \prod_{k=j}^{n} N_{\alpha_k,\alpha_{\max}-\alpha_1-\dots-\alpha_k,1} \xi_k \right)\\
				&\qquad\qquad\cdot
					\cdot x_{\alpha_1+\dots+\alpha_n}(N_{\alpha_1+\dots+\alpha_{n-1},\alpha_n,1} \zeta_{n-1} \xi_n)\\
				&\qquad\qquad\cdot
					x_{\alpha_{\max}}(N_{\alpha_1+\dots+\alpha_{n-1},\alpha_n,2} \zeta_{n-1}^2 \xi_n)).
			\end{aligned}		
		\end{equation}
		Now note that the element $x_{\alpha_{\max}-\alpha_1-\dots-\alpha_{j-1}}(*)$, $j \leq i-1$
		centralizes the set $\cup_{k \geq i} X_{\alpha_1+\dots+\alpha_k} \cdot X_{\alpha_{\max}-\alpha_1-\dots-\alpha_{i-1}}$.
		Therefore, we may rewrite \eqref{eq:lemma:ai:skin:max:big''} as
		\begin{equation}
			\label{eq:lemma:ai:skin:max:big'''}
			\begin{aligned}\relax
				[b,c] 
				&= 
				\left(\prod_{i=2}^{n-1}
					x_{\alpha_1+\dots+\alpha_i}(N_{\alpha_1+\dots+\alpha_{i-1},\alpha_i} \zeta_{i-1} \xi_i)\right.\\
				&\qquad\qquad\cdot
					x_{\alpha_{\max}-\alpha_1-\dots-\alpha_{i-1}}
						(N_{\alpha_{\max}-\alpha_1-\dots-\alpha_i,\alpha_i} \eta_i \xi_i\\
				&\qquad\qquad\qquad
						- \sum_{j=i+1}^{n-1} \eta_j \prod_{k=i}^j N_{\alpha_k,\alpha_{\max}-\alpha_1-\dots-\alpha_k} \xi_k\\
				&\qquad\qquad\qquad
						- \zeta_{n-1} \cdot \prod_{k=i}^{n} N_{\alpha_k,\alpha_{\max}-\alpha_1-\dots-\alpha_k,1} \xi_k)
						\left.\vphantom{\prod_{i=2}^{n-1}}\right) \\
				&\qquad\qquad\cdot
					x_{\alpha_1+\dots+\alpha_n}(N_{\alpha_1+\dots+\alpha_{n-1},\alpha_n,1} \zeta_{n-1} \xi_n)\\
				&\qquad\qquad\cdot
					x_{\alpha_{\max}}(N_{\alpha_{\max}-\alpha_1,\alpha_1} \eta_1 \xi_1 + 
						N_{\alpha_1+\dots+\alpha_{n-1},\alpha_n,2} \zeta_{n-1}^2 \xi_n\\
				&\qquad\qquad\qquad
						-\sum_{j=2}^{n-1} \eta_j \prod_{k=1}^j N_{\alpha_k, \alpha_{\max}-\alpha_1-\dots-\alpha_k} \xi_k\\
				&\qquad\qquad\qquad
						- \zeta_{n-1} \prod_{k=1}^n N_{\alpha_k, \alpha_{\max}-\alpha_1-\dots-\alpha_k,1} \xi_k).
			\end{aligned}		
		\end{equation}		

		Comparing \eqref{eq:lemma:ai:skin:max:1} with \eqref{eq:lemma:ai:skin:max:big'''} we get a system of equations
		which always has a solution. We will construct it in three steps. First, considering the parameters
		of $x_{\alpha_1+\dots+\alpha_n}(*)$ we see that
		\begin{align*}
			N_{\alpha_1+\dots+\alpha_{n-1},\alpha_n,1} \zeta_{n-1} \xi_n = u_n
		\end{align*}
		Set $\zeta_{n-1}=1$ and $\xi_n = (N_{\alpha_1+\dots+\alpha_{n-1},\alpha_n,1})^{-1} u_n$.
		Next, we continue by induction on from $i=n-1$ downwards to $i=2$. On each step
		we have the variables $\xi_{i+1},\dots,\xi_n,$ $\zeta_i,\dots,\zeta_{n-1},$ and $\eta_{i+1},\dots,\eta_{n-1}$
		already set. Considering the parameters of $x_{\alpha_1+\dots+\alpha_i}(*)$ and
		$x_{\alpha_{\max}-\alpha_1-\dots-\alpha_{i-1}}(*)$ we get the equations
		\begin{align*}
			N_{\alpha_1+\dots+\alpha_{i-1},\alpha_i} \zeta_{i-1} \xi_i &= u_i \\
			N_{\alpha_{\max}-\alpha_1-\dots-\alpha_i,\alpha_i} \eta_i \xi_i
			- \xi_i P_i(\xi_{i+1},\dots,\xi_n, \eta_{i+1},\dots,\eta_{n-1},\zeta_{n-1})
			&=v_i,
		\end{align*}
		where $P_i$ is some polynomial. Set $\xi_i = 1$, $\zeta_{i-1} = u_i (N_{\alpha_1+\dots+\alpha_{i-1},\alpha_i})^{-1}$
		and 
		\[
			\eta_i = (v_i+P_i(\xi_{i+1},\dots,\xi_n, \eta_{i+1},\dots,\eta_{n-1},\zeta_{n-1}))/
				(N_{\alpha_{\max}-\alpha_1-\dots-\alpha_i,\alpha_i}).
		\]
		Up till now, we have the variables $\xi_2,\dots,\xi_n$, $\zeta_1,\dots,\zeta_{n-1}$ and
		$\eta_2,\dots,\eta_{n-1}$ set. The last equation we have to satisfy comes from 
		considering the parameters of $x_{\alpha_{\max}}(*)$:
		\begin{align*}
			N_{\alpha_{\max}-\alpha_1,\alpha_1} \eta_1 \xi_1 				
						N_{\alpha_1+\dots+\alpha_{n-1},\alpha_n,2} \zeta_{n-1}^2 \xi_n
			-\xi_1 P_1(\xi_2,\dots,\xi_n, \eta_2,\dots,\eta_{n-1},\zeta_{n-1}) = v_1,
		\end{align*}
		where, again, $P_1$ is a polynomial. Set $\xi_1=1$ and 
		\[
			\eta_1 = (v_1 - N_{\alpha_1+\dots+\alpha_{n-1},\alpha_n,2} \zeta_{n-1}^2 \xi_n + P_1(\xi_2,\dots,\xi_n, \eta_2,\dots,\eta_{n-1},\zeta_{n-1}))/(N_{\alpha_{\max}-\alpha_1,\alpha_1}).
		\]
		Thus $a = [b,c]$. By Lemma \ref{lemma:ai:only:skin}, $\phi(b) = b \cdot x_{\alpha_{\max}}(*)$.
		Further, by Lemma \ref{lemma:ai:prod:of:simple} $\phi(c)=c x_{\alpha_{\max}}(*)$. Thus,
		\[
			\phi(a) =\phi([b,c]) = [\phi(b),\phi(c)] = [b x_{\alpha_{\max}}(*), c x_{\alpha_{\max}}(*)]
			= [b,c] = a.
		\]
	\end{proof}
\end{lemma}

\begin{theorem}
	\label{th:ai:is:central}
	Assume $\rk \Phi \geq 4$. Let $\phi$ be an almost identity $\PC$-map. Then $\phi$ is a central
	$\PC$-map.
	\begin{proof}
		Pick a matrix $a \in \Up$. By Lemma \ref{lemma:ai:giblets} we know that 
		$a_{ij} = \phi(a)_{ij}$ for all $2 \leq i<j \leq 2n-1$.
		Now fix some index $j \neq 1,2n$ and consider the commutator
		\[
			b = [a, x_{\alpha_{\max}-\alpha_{1j}}(1)] \in [\Up, \U_1] \le \U_1^{(2)}.
		\]
		Note that $\alpha_{\max}-\alpha_{1,j} = \alpha_{1,-j}$.
		By direct calculation,
		\begin{align*}
			b = [a, x_{\alpha_{\max}-\alpha}(1)] 
			&= 
			(e + a_{*1} a'_{2n+1-j,*} \pm a_{*j} a'_{2n,*}) x_{\alpha_{\max}-\alpha_{1j}}(-1) \\
			&= 
			(e + e_{*1} a'_{2n+1-j,*} \pm a_{*j} e_{2n,*}) x_{\alpha_{\max}-\alpha_{1j}}(-1).
		\end{align*}
		In particular, $b_{1,2n} = \pm 2 a_{1j} \pm a_{2n+1-j,j}$
		Considering $\phi(b)_{1,2n} = [\phi(a), x_{\alpha_{\max}-\alpha_{1j}}(1)]$ we get in the same way 
		$\phi(b)=\pm 2 \phi(a)_{1j} \pm \phi(a)_{2n+1-j,j}$. By Lemma
		\ref{lemma:ai:skin:max} we have $\phi(b)=b$ and by Lemma \ref{lemma:ai:giblets}
		$a_{2n+1-j,j} = \phi(a)_{2n+1-j,j}$. Thus $a_{1j} = \phi(a)_{1j}$ for all $2 \leq j \leq 2n-1$.
		It follows form the fact that $a$ and $\phi(a)$ are symplectic that $a_{j,2n}=\phi(a)_{j,2n}$ for
		all $2 \leq j \leq 2n-1$. 	
		Summing up, $a = \phi(a) x_{\alpha_{\max}}(*)$.
	\end{proof}
\end{theorem}

\noindent{\textbf{Proof of Theorem \ref{th:main}.}}
	By Theorem \ref{theorem:up:to:ai}, $\phi$ can be decomposed as follows:
	\[
		\phi = \Int_C \circ E \circ Int_{C'} \circ D \circ Q \circ Z \circ \tau \circ A,
	\]
	where $\Int_C$ and $\Int_{C'}$ are inner, $E$ an extremal, $D$ a diagonal, $Q$ a quasi-diagonal,
	$\tau$ a field automorphisms; $Z$ -- a central $\PC$-map and $A$ -- an almost identity $\PC$-map.
	By Theorem \ref{th:ai:is:central}, $A$ is central.
	\qed
	
\bibliography{../references-eng}

\end{document}